\documentclass{svjour3}       
\usepackage{amsmath,amssymb}
\usepackage{float}
\usepackage{verbatim,color}
\usepackage{float}
\usepackage{hyperref}
\RequirePackage{fix-cm}
\usepackage{xcolor}
\newcommand{\Allen}[1]{{{\color{black} {#1}}}}
\newcommand{\deaa}[1]{{{\color{black} {#1}}}}
\newcommand{\Alldeaa}[1]{{{\color{black} {#1}}}}
\newcommand{\FinalEdits}[1]{{{\color{black} {#1}}}}
\smartqed  
\usepackage{graphicx}
\begin{document}

\title{A Hermite Method with a Discontinuity Sensor for Hamilton-Jacobi Equations\thanks{AAL was supported by NSF Grant DGE-1650115. DA was supported, in part, by NSF Grant DMS-1913076. Any opinions, findings, and conclusions or recommendations expressed in this material are those of the authors and do not necessarily reflect the views of the National Science Foundation}
}


\author{Allen Alvarez Loya         \and
        Daniel Appel\"{o} 
}


\institute{Allen Alvarez Loya \at
              Department of Applied Mathematics, University of Colorado at Boulder\\
              Engineering Center, ECOT 225\\
              526 UCB\\
              Boulder, CO 80309-0526\\
              Tel.: 714-864-0760\\
              \email{allen.alvarezloya@colorado.edu}           
           \and
Daniel Appel\"{o} \at Department of Computational Mathematics, Science \& Engineering and Department of Mathematics, Michigan State University, East Lansing, USA.  \\
\email{appeloda@msu.edu}
}

\date{Received: date / Accepted: date}

 \bibliographystyle{plain}
\maketitle

\begin{abstract}
We present a Hermite interpolation based partial differential equation solver for
Hamilton-Jacobi equations. Many Hamilton-Jacobi equations have a nonlinear
dependency on the gradient, which gives rise to discontinuities in the
derivatives of the solution, resulting in kinks. We built our solver with two goals in mind: 1) high order accuracy in smooth regions and 2) sharp resolution of
kinks. To achieve this, we use Hermite interpolation with a smoothness sensor. The
degrees-of-freedom of Hermite methods are tensor-product Taylor polynomials of
degree $m$ in each coordinate direction. The method uses $(m+1)^d$ degrees of
freedom per node in $d$-dimensions and achieves an order of accuracy $(2m+1)$
when the solution is smooth. To obtain sharp resolution of kinks, we sense the
smoothness of the solution on each cell at each timestep. If the solution is smooth, we
march the interpolant forward in time with no modifications. When our method
encounters a cell over which the solution is not smooth, it introduces artificial viscosity
locally while proceeding normally in smooth regions. We show through numerical
experiments that the solver sharply captures kinks once the solution losses continuity
in the derivative while achieving $2m+1$ order accuracy in smooth regions.  
\end{abstract}


\section{Introduction}  
We consider the time-dependent Hamilton-Jacobi (HJ) equation 
\begin{equation}
\varphi_t + H(\varphi_{x_1},\varphi_{x_2},\dots,\varphi_{x_d}) = 0, \quad \varphi(\mathbf{x},0) = \varphi^0(\mathbf{x}), \quad \mathbf{x} \in \Omega \in \mathbb{R}^d,
\label{HJ}
\end{equation}
with periodic boundary conditions on $\partial \Omega$. HJ equations appear in many applications, e.g., optimal control, differential games, image processing and the calculus of variations. The solution of an HJ equation may develop a discontinuity in the derivative even when the initial data is smooth. As in conservation laws, the unique solution can be singled out by the use of  viscosity solutions \cite{crandall1983viscosity,crandall1984some}. In particular, the viscosity solution gives requirements on the solution at points of discontinuity that allow us to find the unique physically relevant solution. Design of methods for HJ equations that: a) converge to the viscosity solution as the grid is refined even in the presence of kinks, and b) maintain high-order accuracy in smooth regions and here we only mention a few.
The first methods to try to accomplish this were introduced in \cite{souganidis1985approximation,crandall1984two}. Essentially non-oscillatory (ENO) \cite{osher1991high} or weighted ENO (WENO) \cite{jiang2000weighted,zhang2003high} have been developed to solve the HJ equation. These methods have been shown to work efficiently, but require a large stencil size in order to obtain high-order accuracy. More recently, there has been work done by using discontinuous Galerkin methods to solve conservation laws. This idea was introduced in \cite{osher1991high}, where the authors exploited the fact that the gradient of the solution satisfies a conservation law system, and used a standard discontinuous Galerkin method to advance the solution in time and then recovered the solution from the derivatives. In reference \cite{cheng2007discontinuous} the authors developed a discontinuous Galerkin method for directly solving \eqref{HJ}, when the Hamiltonian is linear or convex, this eliminates the need to solve systems in the multidimensional case. This work was later improved upon in \cite{cheng2014new}. The improvement extends the method to approximate solutions to \eqref{HJ} when the Hamiltonian is not convex. Other improvements include avoiding reconstruction of the solution across elements by utilizing the Roe speed at the cell interfaces and adding an entropy fix inspired by Harten and Hyman \cite{harten1983self}. 

In \cite{qiu2005hermite} Qiu and Shu use WENO methods together with Hermite interpolating polynomials (HWENO). These methods are derived from the original WENO methods, but both the function and the first derivative are evolved and used in the reconstruction of the solution. An important advantage of HWENO over WENO is that a more compact stencil may be used for the same order of accuracy. There have been several expansions on this work \cite{qiu2007hermite,zheng2017finite,tao2015high,zahran2016seventh}, where in \cite{zahran2016seventh} the authors develop a seventh order method, which is the same order of accuracy we obtain using three derivatives. The success of HWENO methods has inspired us to build a Hermite method solver for \eqref{HJ}. Using Hermite methods, we may achieve an arbitrary order while keeping a compact stencil. 

Hermite methods were first studied in \cite{goodrich2006hermite}, where the authors use Hermite methods to solve hyperbolic initial-boundary value problems. Stability of the method and convergence was proved and various numerical examples were provided. Several adaptations of the original Hermite method have been developed. In \cite{appelo2018hermite} the authors use Hermite interpolants to solve the wave equation using dissipative and conservative formulations. A hybrid Hermite-discontinuous Galerkin method was used in \cite{chen2014hybrid}, where the authors approximated the solution of Maxwell's equations. In \cite{cons_hermite} Hermite methods for hyperbolic conservation laws were considered, where the entropy viscosity by Guermond et al., \cite{Guermond2008801}, was adopted to the Hermite framework.  The current paper presents the first Hermite method for Hamilton-Jacobi equations. 

The challenge in solving Hamilton-Jacobi equations stems from the nonlinear dependency on the gradient, which gives rise to discontinuities in the derivatives of the solution, resulting in kinks. To resolve the kinks, we sense the smoothness of the solution on each cell at each timestep by adopting the sensor introduced by Persson and Peraire in  \cite{persson2006sub} and later refined by Klockner, Warburton and Hesthaven in \cite{klockner2011viscous}. Following \cite{persson2006sub,klockner2011viscous} we make no modifications  if the solution is smooth but when the solution fails to be smooth we locally introduce artificial viscosity. Through numerical experiments we demonstrate that the solver has the following properties: 1) high order accuracy in smooth regions and 2) sharp resolution of kinks. \deaa{Although (as most high order methods) we cannot expect the order of accuracy to be $(2m+1)$ when the solution becomes non-smooth we observed in \cite{AppHaglinear} that the constant in from the the $h$ term in the error expansion can be significantly smaller when the formal order is high.}
 
The rest of this paper is organized as follows: in Section 2 we introduce the numerical scheme for the one-dimensional and two-dimensional HJ equations. Section 3 explains how we sense the smoothness of each element. Section 4 is devoted to the discussion of numerical experiments. \FinalEdits{Code for the numerical examples can be found in the github repository \url{https://github.com/allenalvarezloya/Hermite_HJ}}. Conclusions and future work are discussed in Section 5.

\section{Hermite Methods}
A Hermite method of order $2m+1$ approximates the solution to a PDE by element-wise polynomials that interpolate the solution and derivatives up to degree $m$ at the element interfaces. The time evolution is done locally from the center of the element. In Hermite methods, the degrees of freedom are the function and the spatial derivatives, or equivalently the Taylor coefficients. The evolution of the polynomials depends on the nature of the PDE.

\subsection{Hermite Method in One Dimension}
We describe a Hermite PDE method in one-dimension for HJ equations. This method closely follows the Hermite solver for conservation laws described in \cite{hagstrom2015solving}. Note that we use this method, with no modifications, while the solution is smooth. Consider the Hamilton-Jacobi equation, with periodic boundary conditions 
\begin{equation}
u_t = \begin{cases} 
\varphi_t +H(\varphi_x) = 0,\\
\varphi(x,0)=\varphi^0(x).\\
\end{cases} 
\label{PDE}
\end{equation}
\subsubsection{Initialization}
We discretize the spatial coordinate into a primal grid and a dual grid
\begin{equation}
x_i = x_l + i h_x, \quad h_x = (x_r - x_l)/n_x,
\label{Grid}
\end{equation}
where $i = 0,\dots, n_x$ for the primal grid and $i = 1/2, \dots , n_x - 1/2$ for the dual grid. The first step initializes the degrees of freedom by setting them to be the scaled derivatives at each primal grid point. That is, 
\[ c_{l,i} = \left.\frac{h_x^l}{l!} \frac{d^l u(x,0)}{dx^l} \right|_{x=x_i}. \]
We then interpolate scaled derivative data to obtain the piecewise degree-$2m+1$ polynomial
\begin{equation}
v_{i+\frac{1}{2}}(x,0) = \sum \limits_{l = 0}^{2m+1} d_{l,i+\frac{1}{2}} \left(\frac{x - x_{i+\frac{1}{2}}}{h_x }\right)^l,
\end{equation}
where the coefficients $d_{l,i}$ are uniquely determined by the interpolation conditions
\begin{equation}
\frac{h_x^l}{l!}\frac{\partial^l }{\partial x^l} v_{i+\frac{1}{2}}(x_i) = c_{l,i}, \quad \frac{h_x^l}{l!}\frac{\partial^l }{\partial x^l} v_{i+\frac{1}{2}}(x_{i+1}) = c_{l,i+1},  
\end{equation}
for $l = 0, \dots, m$.
\subsubsection{Evolution}
We treat the coefficients in the Hermite interpolant as functions of time. That is, 
\begin{equation}
v_{i+\frac{1}{2}}(x,t) = \sum \limits_{l = 0}^{2m+1} d_{l,i+\frac{1}{2}}(t) \left(\frac{x - x_{i+\frac{1}{2}}}{h_x }\right)^l.
\label{PDE_Hermite_expansion}
\end{equation}
For our PDE $\varphi_t = -H(\varphi_x)$ we substitute \eqref{PDE_Hermite_expansion}:
\begin{equation}
\frac{\partial v_{i+\frac{1}{2}}(x,t) }{\partial t } =  \sum \limits_{l = 0}^{2m+1} d_{l,i+\frac{1}{2}}^{'}(t) \left(\frac{x - x_{i+\frac{1}{2}}}{h_x }\right)^l = -H(v_x(x,t)).
\label{PDE_Hermite_expansion_substitution}
\end{equation}
We can differentiate \eqref{PDE_Hermite_expansion_substitution} $k$ times in space and evaluate at $x = x_{i+\frac{1}{2}}$ to obtain
\begin{equation}
\frac{k!}{h_x^k} d^{'}_{k,i+\frac{1}{2}} (t) = -\left. \frac{\partial^k}{\partial x^k} H(v_x(x,t)) \right |_{x = x_{i + \frac{1}{2}}}.
\label{DiffOfHamiltonian}
\end{equation}
When $H$ is nonlinear the differentiation of the RHS can spawn new terms by the product rule. We avoid this by approximating the RHS by a Taylor polynomial of degree $2m+1$ 
\begin{equation}
H(v_x) \approx \sum \limits_{l=0}^{2m+1} b_{l,i+\frac{1}{2}}(t) \left (\frac{x - x_{i + \frac{1}{2}}}{h_x} \right )^l.
\label{Taylor_expansion}
\end{equation}
for which the differentiation is straight forward. With this approximation we carry out the differentiation in \eqref{DiffOfHamiltonian} and obtain the local system of ODEs
\begin{equation}
d^{'}_{k,i+\frac{1}{2}}(t) = b_{k,i+\frac{1}{2}}(t), \quad k = 0, \dots , 2m+1.
\label{ODEs1D}
\end{equation}
We can solve this system to evolve our approximate solution using a standard one step ODE solver as described in \cite{HagApp07}. Before we evolve we must find the Taylor coefficients $b_{k,i+\frac{1}{2}}(t).$ The computation of the coefficients, $b_{k,i+\frac{1}{2}}$, is problem specific and depends on the form of the Hamiltonian, $H$.

\subsubsection{Polynomial Approximation of Non-linear Hamiltonians}
\deaa{The method described above relies on the computation of the Taylor series coefficients of functions of polynomials. It is convenient and efficient to compute these coefficients through recursions that rely on Cauchy products rather than through the use of the higher order chain rule (the Fa\`{a} di Bruno's formula) which have combinatorial complexity. Our discussion here closely follows \cite{neidinger2013efficient}.}

We describe this process for the PDE that is used in example 3 in the numerical experiments section below
\begin{align*} 
\varphi_t = \cos(\varphi_x + 1).
\end{align*}

A Taylor series expansion for the right hand side of this equation can be obtained by \Allen{the chain rule}. If in general the functions $f(x), w(x)$ and $u(x)$ have Taylor series expansions
\begin{align*}
f(x) &= \sum \limits_{k=0}^{\infty} F_k \left(\frac{x - x_i}{h} \right)^k,\\
w(x) &= \sum \limits_{k=0}^{\infty} W_k \left(\frac{x - x_i}{h} \right)^k,\\
u(x) &= \sum \limits_{k=0}^{\infty} U_k \left(\frac{x - x_i}{h} \right)^k,
\end{align*}
then for non-linearities that satisfy a differential equation $f'(x) = w(x)u'(x)$ we can compute $F_k, k = 1,2,\dots $ using the formula
\begin{align*}
F_k = W_0 U_k + \frac{1}{k} \sum \limits_{j = 0}^{k-1} j U_j W_{k-j}.
\end{align*}

Here we are interested in approximating $\cos(\varphi_x+1)$ and thus note that the functions $s(x) = \sin(u(x))$ and $c(x) = \cos(u(x))$ satisfy $s'(x) = c(x)u'(x)$ and $c'(x) = s(x)u'(x)$ simultaneously. Thus, we use both relationships and the formula above to obtain $\cos(\varphi_x+1)$. Precisely, first we compute 
\[ (\varphi_x)_k = \left\{ \begin{array}{ll}
         \frac{k+1}{h} \varphi_k & \mbox{if $k = 0 \dots 2m$},\\
        0 & \mbox{if $k = 2m+1$},\end{array} \right. \] 
followed by
\[ \Allen{U_k} = \left\{ \begin{array}{ll}
         (\varphi_x)_k + 1& \mbox{if $k = 0$},\\
        (\varphi_x)_k & \mbox{if $k = 1,\dots, 2m+1$}.\end{array} \right. \] 
\Allen{Note that the coefficients, $U_k$, are simply computed from the Taylor expansion of $\varphi$.} We then set $C_0 = \cos(U_0)$ and $S_0 = \sin(U_0)$ and update the remaining coefficients using the recursion 
\begin{align*}
S_k &= C_0U_k+\texttt{ddot}(C,U,k),\\
C_k &= -S_0U_k-\texttt{ddot}(S,U,k),
\end{align*}
where  $u = \Allen{[U_0, U_1, \dots, U_{2m+1}]^T}$ , $c = \Allen{[C_0, C_1, \dots, C_{2m+1}]^T}$,  $s = \Allen{[S_0, S_1, \dots, S_{2m+1}]^T}$ and \verb+ddot+ is the function given by
\begin{align*}
\texttt{ddot}(A,B,k) = \frac{1}{k} \sum \limits_{j = 0}^{k-1} j A_j B_{k-j}.
\end{align*}
The $C_k$ are the Taylor coefficients of the polynomial approximating $\cos(\varphi_x+1)$.

\begin{figure}[htb]
\begin{center}
\includegraphics[width=0.49\textwidth]{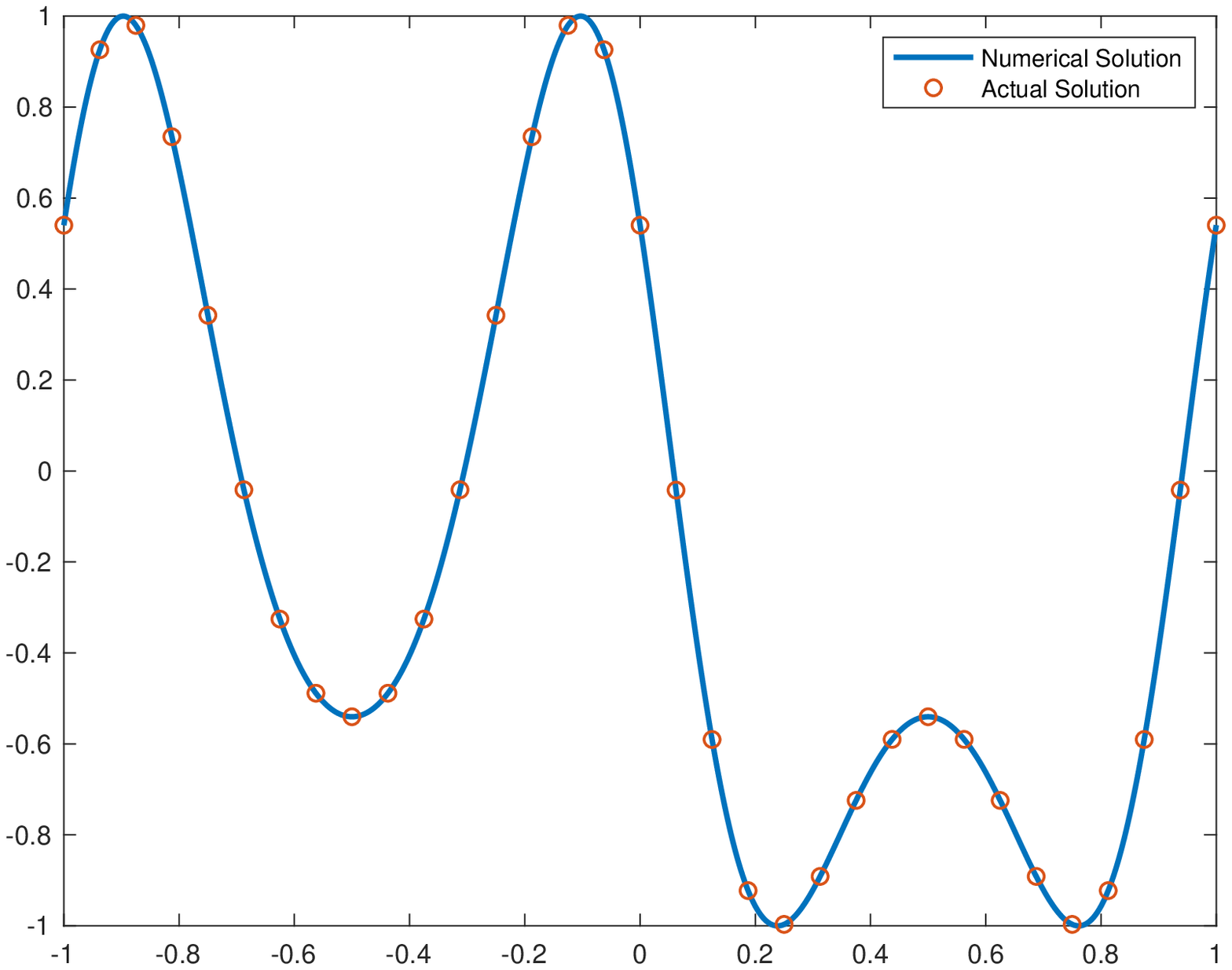} 
\includegraphics[width=0.49\textwidth]{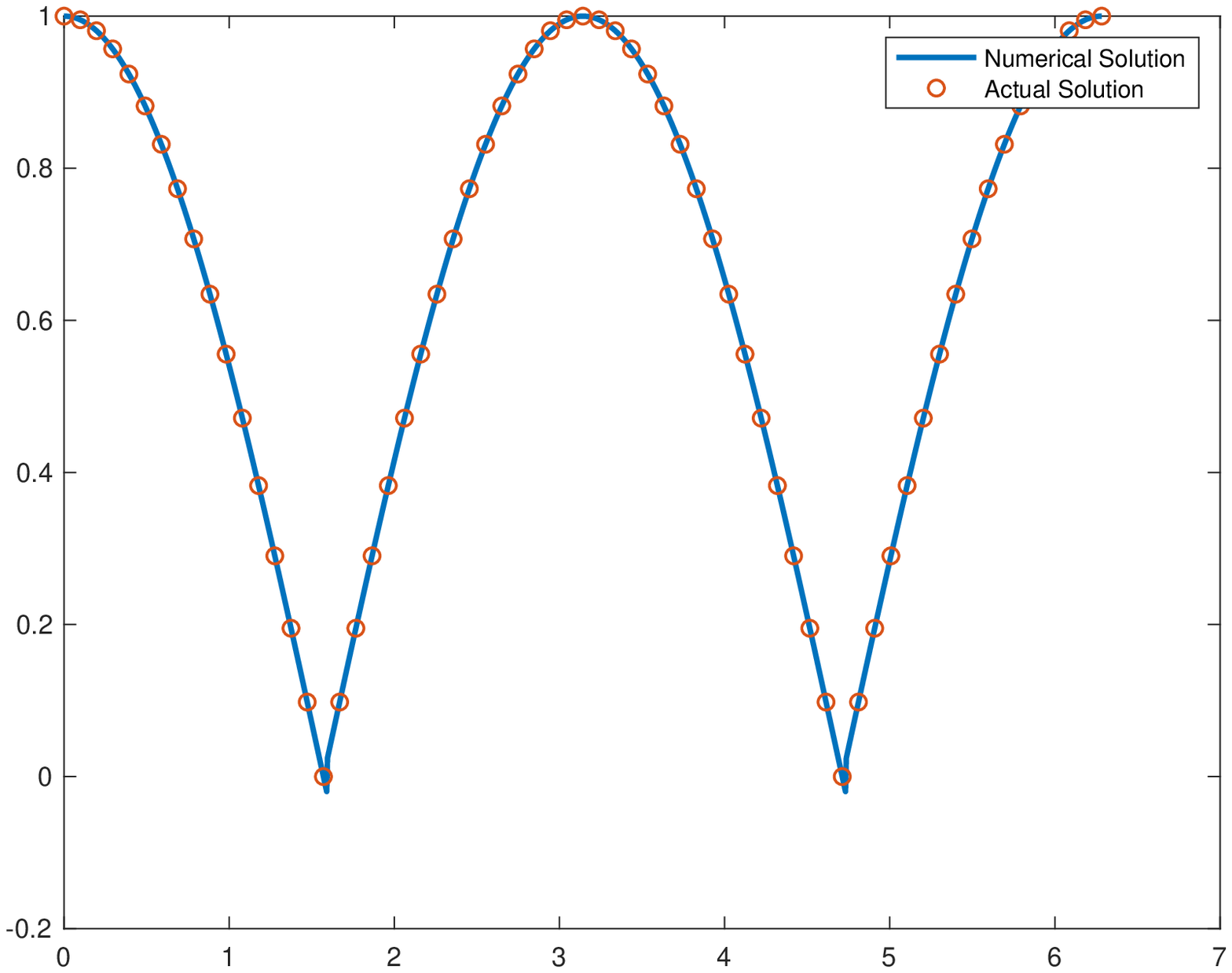} 
\caption{On the left we plot the numerical approximation to $H(\varphi_x) = -\cos(\varphi_x+1)$ as a solid line and sample the actual function as circles. On the right we plot the numerical approximation to $H(\varphi_x) = -|\varphi_x|$ as a solid line and sample the actual function as circles. In both cases we used $m = 3$ derivatives. \label{fig:TaylorExample}}
\end{center}
\end{figure}

In Table \ref{tab:TaylorTable1} we display the $L_1,L_2$ and $L_{\infty}$ norm errors for $H(\varphi_x) = -\cos(\varphi_x+1)$ using $x_l = 0$, $x_r = 2\pi$ using $n_x = 20, 40, 80 $ and 160 gridpoints. As can be seen from the estimated rates of convergence the procedure is as accurate as the Hermite interpolation ($h^{2m+2}$). To the left in Figure \ref{fig:TaylorExample} we plot the Taylor series approximation which can be seen to be accurate.  
\begin{table}[ht]
\caption{Errors in approximating $H(\varphi_x) = -\cos(\varphi_x+1)$ using Taylor series approximation recursions in the $L_1$, $L_2$ and $L_{\infty}$ norms at time are displayed along with estimated rates of convergence. Note that we obtain the expected rate of $2m+2$ where $m$ is the number of derivatives used in approximation.\label{tab:TaylorTable1}} 
\begin{center} 
\begin{tabular}{c c c c c c c} 
\hline 
n & $L_1$ error & Convergence & $L_2$ error & Convergence & $L_{\infty}$ error & Convergence \\ 
\hline 
& & & $m = 2$ & & &\\
\hline 
20 & 1.18e-05 & - & 1.95e-05 & - & 5.05e-05 & - \\ 
40 & 1.75e-07 & 6.08 & 3.07e-07 & 5.99 & 9.89e-07 & 5.68 \\ 
80 & 2.74e-09 & 5.99 & 4.86e-09 & 5.98 & 1.61e-08 & 5.94 \\ 
160 & 4.44e-11 & 5.95 & 7.94e-11 & 5.94 & 2.71e-10 & 5.90 \\ 
\hline 
& & & $m = 3$ & & &\\
\hline 
20 & 1.42e-07 & - & 2.64e-07 & - & 7.65e-07 & - \\ 
40 & 5.69e-10 & 7.96 & 1.04e-09 & 7.99 & 3.37e-09 & 7.83 \\ 
80 & 2.19e-12 & 8.02 & 4.06e-12 & 8.00 & 1.32e-11 & 8.00 \\ 
160 & 8.93e-15 & 7.94 & 1.59e-14 & 8.00 & 5.20e-14 & 7.99 \\ 
\hline 
\end{tabular} 
\end{center} 

\end{table} 

In the case of a Hamiltonian with an absolute value function the above procedure cannot be used. We describe how to handle such cases using the Eikonal equation
\begin{align*}
\varphi_t + |\varphi_x| = 0,
\end{align*}
as an example. Here, we need the expansion of $|\varphi_x|$. If $\varphi_x$ is not sign definite, then we can not use the method above to obtain a Taylor series expansion and instead we proceed by computing
\[ (\varphi_x)_k = \left\{ \begin{array}{ll}
         \frac{k+1}{h} \varphi_k & \mbox{if $k = 0 \dots 2m$},\\
        0 & \mbox{if $k = 2m+1$}.\end{array} \right. 
\] 
At the cell center the value of the function is given by the leading coefficient. That is, $\varphi_x(x_i) = \varphi_0$. Therefore, as an approximation we take 
\[ 
|(\varphi_x)_k| = \left\{ \begin{array}{ll}
         (\varphi_x)_k & \mbox{if $\varphi_0 \geq 0$},\\
        -(\varphi_x)_k & \mbox{if $\varphi_0 < 0$}.\end{array} \right. 
\] 

In Table \ref{tab:TaylorTable1} we display the $L_1,L_2$ and $L_{\infty}$ norm errors for $H(\varphi_x) = |\varphi_x|$ using $x_l = 0$, $x_r = 2\pi$ with $n_x = 20, 40, 80 $ and 160 gridpoints. We show the estimated rates of convergence for each norm. Of course here we cannot expect the full order of convergence as the nonlinearity is not smooth. To the right in Figure \ref{fig:TaylorExample} we plot the Taylor series approximation, which can be seen to be an accurate approximation to the true Hamiltonian. \Allen{Since the degrees of freedom are used from a single point to evaluate the solution on the entire cell a change in the sign of $\varphi_x$ results in an error of either $2\varphi_x$ or -$2\varphi_x$ over the interval where the sign change occurs to the end of the cell. This implies that the error is $\mathcal{O}(h)$ over the cells where a change in sign of $\varphi_x$ occurs.}  
\begin{table}[ht]
 \caption{Errors in approximating $H(\varphi_x) = |\varphi_x|$ using Taylor series in the $L_1$, $L_2$ and $L_{\infty}$ norms at time are displayed along with estimated rates of convergence. Note that the function we are trying to interpolate is nonlinear, thus we do not obtain the expected rate of $2m+2$ where $m$ is the number of derivatives used in approximation.} 
\begin{center} 
\Allen{
 \begin{tabular}{c c c c c c c } 
\hline 
 n & $L_1$ error & Convergence & $L_2$ error & Convergence & $L_{\infty}$ error & Convergence \\ 
 \hline 
& & & $m = 2$ & & &\\
\hline 
 20 & 6.15e-02 & - & 1.20e-01 & - & 3.13e-01 & - \\ 
40 & 1.54e-02 & 2.00 & 4.26e-02 & 1.50 & 1.57e-01 & 1.00 \\ 
80 & 3.85e-03 & 2.00 & 1.51e-02 & 1.50 & 7.85e-02 & 1.00 \\ 
160 & 9.64e-04 & 2.00 & 5.33e-03 & 1.50 & 3.93e-02 & 1.00 \\ 
\hline
& & & $m = 3$ & & &\\
\hline
 20 & 6.15e-02 & - & 1.20e-01 & - & 3.13e-01 & - \\ 
40 & 1.54e-02 & 2.00 & 4.26e-02 & 1.50 & 1.57e-01 & 1.00 \\ 
80 & 3.85e-03 & 2.00 & 1.51e-02 & 1.50 & 7.85e-02 & 1.00 \\ 
160 & 9.64e-04 & 2.00 & 5.33e-03 & 1.50 & 3.93e-02 & 1.00 \\ 
\hline 
\end{tabular} 
}
\end{center} 
\label{tab:TaylorTable2}
\end{table} 

\subsubsection{\Allen{Polynomial Approximation of Non-linear Hamiltonians Higher-Dimensions}}
\begin{figure}[hbt]
\begin{center}
\includegraphics[width=0.32\textwidth]{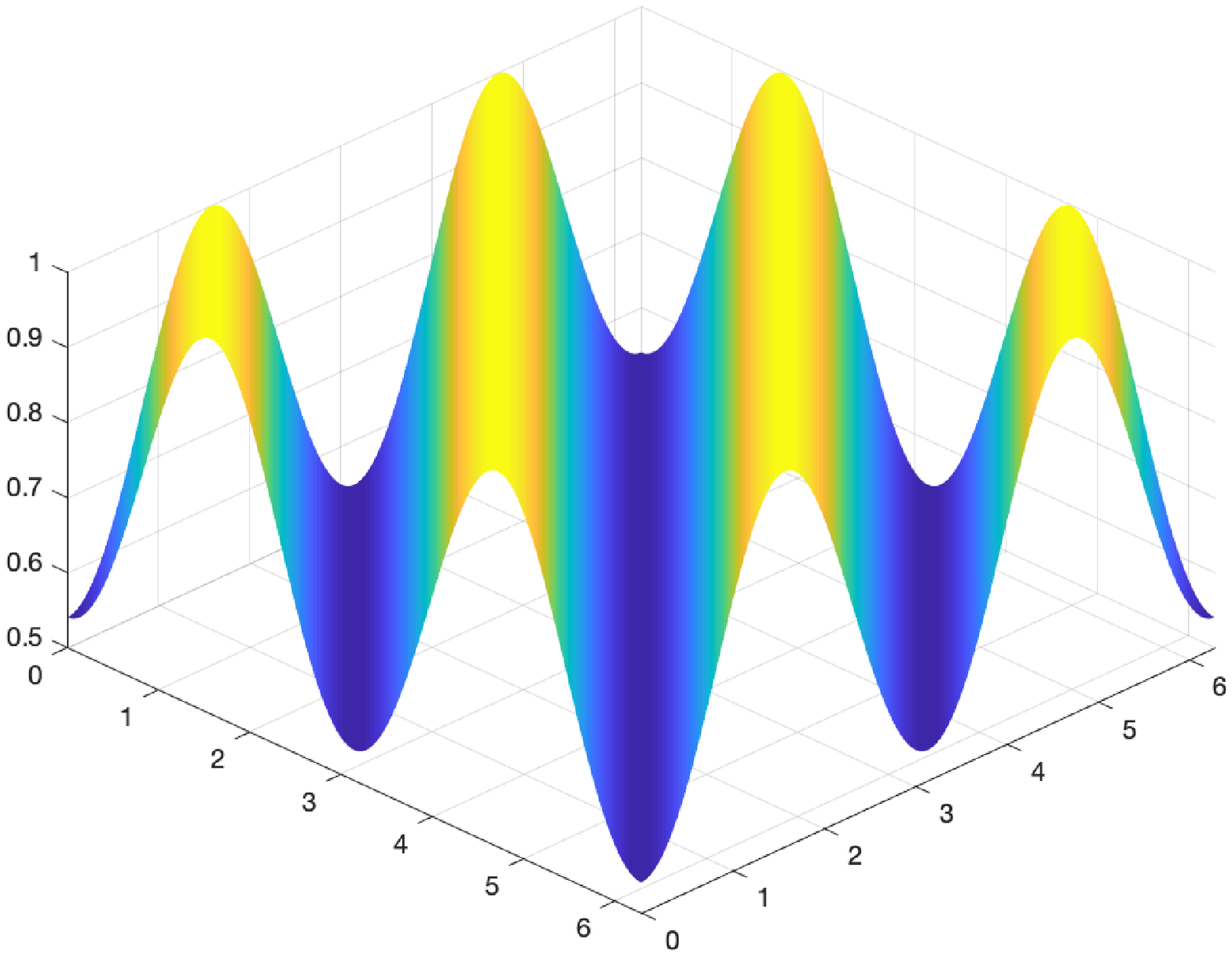} 
\includegraphics[width=0.32\textwidth]{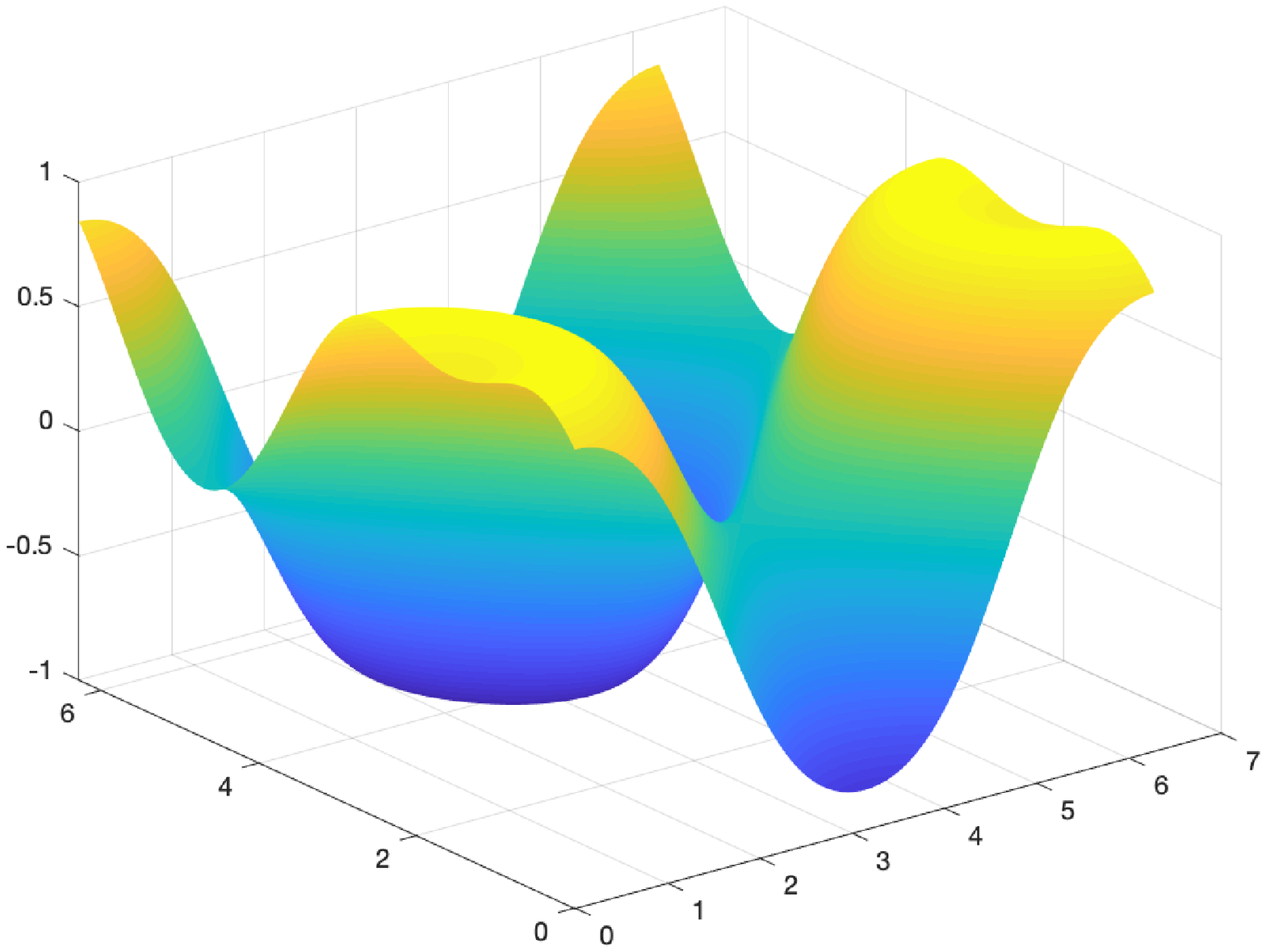}
\includegraphics[width=0.32\textwidth]{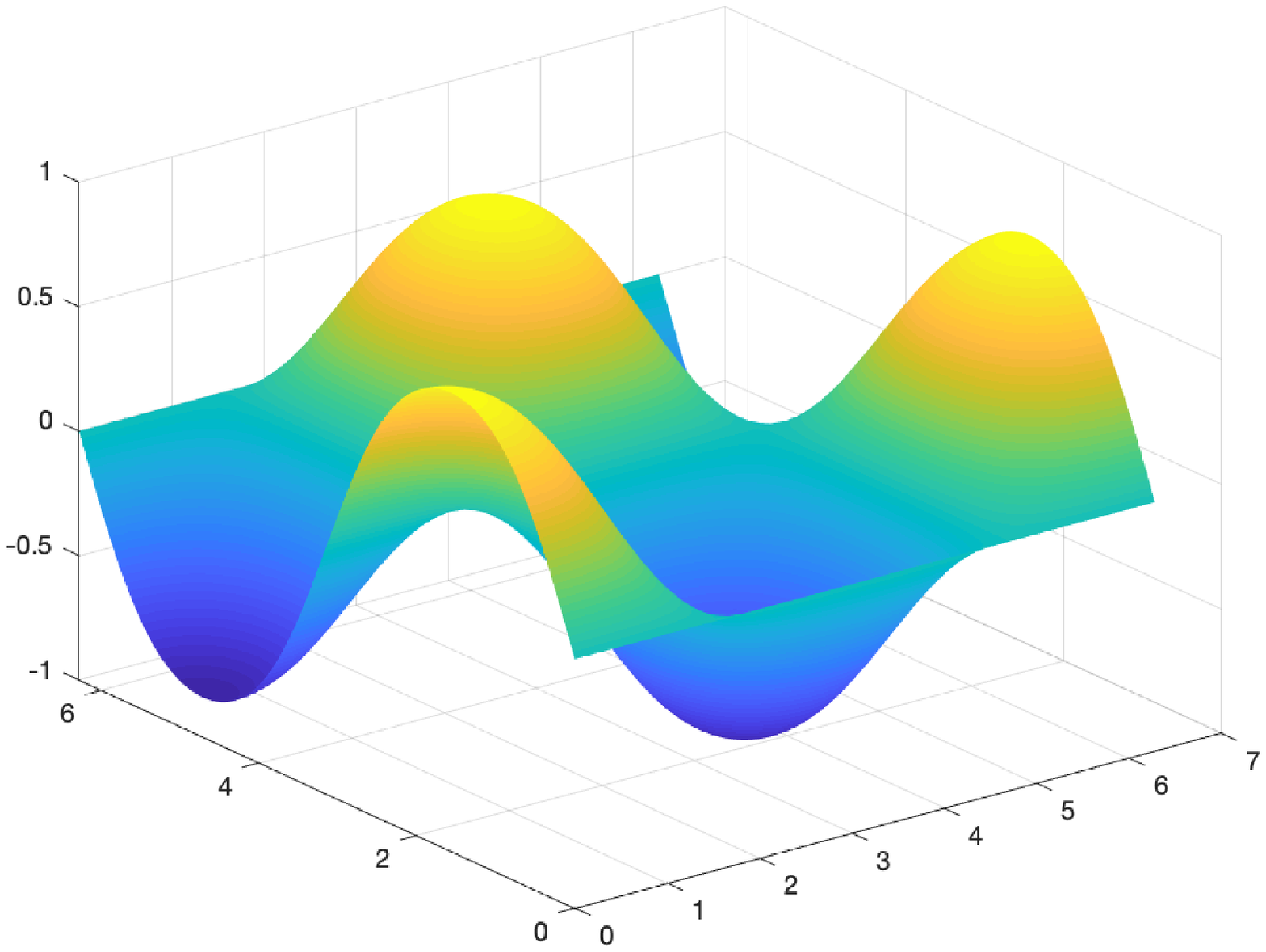} 
\caption{\Allen{We plot the approximations obtained by the Taylor series recursion examples. Left: $H(x,y) = \cos(\cos(x+y))$. Middle: $H(x,y) = \sin(\sin(x) + \cos(y))$. Right: $H(x,y) = \sin(\sin(x)\cos(y))$. For each example $n = 20$ cells were used with $m = 3$ derivatives.}} \label{fig:Taylor2D}
\end{center}
\end{figure}

\Allen{
The computation of the multivariate Taylor coefficients of functions of polynomials is similar to the univariate case. For a composition $H(\mathbf{x}) = f(u(\mathbf{x}))$ where $f:\mathbb{R} \to \mathbb{R}$ is a standard function we have 
\begin{align*}
\frac{\partial h}{\partial x_i}(\mathbf{x}) & = f'(u(\mathbf{x}))\frac{\partial u}{\partial x_i}(\mathbf{x})\\
& = w(\mathbf{x})\frac{\partial u}{\partial x_i}(\mathbf{x}),
\end{align*}
which is similar to the relationship $f'(x) = w(x)u'(x)$ in one-dimension. Using this relationship one can derive the recursion 
\begin{align*}
H_{\mathbf{k}} & = W_{\mathbf{0}}U_{\mathbf{k}} + \texttt{ddot}(V,U,\mathbf{k}).
\end{align*}
Here 
\begin{align*}
\texttt{ddot}(V,U,\mathbf{k}) = \frac{1}{k_p} \sum \limits_{\mathbf{e}_p \leq \mathbf{j} \lvertneqq \mathbf{k}} j_p U(\mathbf{j}) W(\mathbf{k} - \mathbf{j}), 
\end{align*}
where the sum is over all multi-indicies with coordinates $j_p = 1,2,\dots,k_p$ and all other $j_i = 0,1,\dots,k_p$ except when $\mathbf{j} = \mathbf{k}$.

We approximate the functions $H(x,y) = \cos(\cos(x+y))$, $H(x,y) = \sin(\sin(x)+\cos(y))$ and $H(x,y) = \sin(\sin(x)\cos(y))$ using the recursion described above. In each example we estimate the rate of convergence by refining the grid by a factor of two starting with $n = 10$ cells in each coordinate direction and ending with $n = 80$. We compute using $m =2$ and $m = 3$ derivatives. 

For each example we observe the optimal rate of convergence reported in Table \ref{tab:Taylor2Dcoscos}, Table \ref{tab:Taylor2DsinAdd} and Table \ref{tab:Taylor2DsinMult} for $H(x,y) = \cos(\cos(x+y))$, $H(x,y) = \sin(\sin(x)+\cos(y))$ and $H(x,y) = \sin(\sin(x)\cos(y))$, respectively. We display the approximations in Figure \ref{fig:Taylor2D}.
} 
\begin{table}[ht]
 \caption{\Allen{Errors in approximating $f(x,y) = \cos(\cos(x+y))$ using Taylor series in the $L_1$, $L_2$ and $L_{\infty}$ norms at time are displayed along with estimated rates of convergence. We obtain the expected rate of $2m+2$ where $m$ is the number of derivatives used in approximation.}} 
\begin{center} 
\Allen{
 \begin{tabular}{c c c c c c c } 
\hline 
 n & $L_1$ error & Convergence & $L_2$ error & Convergence & $L_{\infty}$ error & Convergence \\ 
\hline 
& & & $m = 2$ & & &\\
\hline
 10 & 4.03e-04 & - & 1.18e-04 & - & 4.14e-05 & - \\ 
20 & 6.23e-06 & 6.02 & 1.83e-06 & 6.01 & 8.33e-07 & 5.63 \\ 
40 & 9.69e-08 & 6.01 & 2.86e-08 & 6.00 & 1.32e-08 & 5.98 \\ 
80 & 1.51e-09 & 6.00 & 4.47e-10 & 6.00 & 2.06e-10 & 6.00 \\
\hline 
& & & $m = 3$ & & &\\
\hline
10 & 9.86e-06 & - & 2.89e-06 & - & 2.11e-06 & - \\ 
20 & 3.82e-08 & 8.01 & 1.09e-08 & 8.04 & 6.98e-09 & 8.24 \\ 
40 & 1.49e-10 & 8.01 & 4.27e-11 & 8.00 & 3.14e-11 & 7.80 \\ 
80 & 5.81e-13 & 8.00 & 1.67e-13 & 8.00 & 1.27e-13 & 7.95 \\ 
\end{tabular} 
}
\end{center} 
\label{tab:Taylor2Dcoscos}
\end{table}

\begin{table}[ht]
 \caption{\Allen{Errors in approximating $f(x,y) = \sin(\sin(x)+\cos(y))$ using Taylor series in the $L_1$, $L_2$ and $L_{\infty}$ norms at time are displayed along with estimated rates of convergence. We obtain the expected rate of $2m+2$ where $m$ is the number of derivatives used in approximation.}} 
\begin{center} 
\Allen{
 \begin{tabular}{c c c c c c c} 
\hline 
 n & $L_1$ error & Convergence & $L_2$ error & Convergence & $L_{\infty}$ error & Convergence \\ 
\hline 
& & & $m = 2$ & & &\\
\hline 
10 & 4.33e-04 & - & 1.17e-04 & - & 8.38e-05 & - \\ 
20 & 6.69e-06 & 6.02 & 1.83e-06 & 6.00 & 1.48e-06 & 5.82 \\ 
40 & 1.04e-07 & 6.01 & 2.85e-08 & 6.00 & 2.31e-08 & 6.00 \\ 
80 & 1.62e-09 & 6.01 & 4.46e-10 & 6.00 & 3.60e-10 & 6.01 \\ 
\hline 
& & & $m = 3$ & & &\\
\hline 
 10 & 8.70e-06 & - & 2.45e-06 & - & 1.70e-06 & - \\ 
20 & 3.35e-08 & 8.02 & 9.53e-09 & 8.01 & 7.95e-09 & 7.74 \\ 
40 & 1.30e-10 & 8.01 & 3.72e-11 & 8.00 & 3.12e-11 & 7.99 \\ 
80 & 5.06e-13 & 8.00 & 1.45e-13 & 8.00 & 1.22e-13 & 8.00 \\ 
\hline 
\end{tabular}
} 
\end{center} 
\label{tab:Taylor2DsinAdd}
\end{table}

\begin{table}[ht]
 \caption{\Allen{Errors in approximating $f(x,y) = \sin(\sin(x)\cos(y))$ using Taylor series in the $L_1$, $L_2$ and $L_{\infty}$ norms at time are displayed along with estimated rates of convergence. We obtain the expected rate of $2m+2$ where $m$ is the number of derivatives used in approximation.}} 
\begin{center}
\Allen{ 
 \begin{tabular}{c c c c c c c } 
\hline 
 n & $L_1$ error & Convergence & $L_2$ error & Convergence & $L_{\infty}$ error & Convergence \\ 
\hline 
& & & $m = 2$ & & &\\
\hline 
10 & 2.11e-04 & - & 7.01e-05 & - & 4.76e-05 & - \\ 
20 & 3.22e-06 & 6.03 & 1.10e-06 & 6.00 & 8.39e-07 & 5.83 \\ 
40 & 5.17e-08 & 5.96 & 1.71e-08 & 6.00 & 1.34e-08 & 5.97 \\ 
80 & 8.10e-10 & 6.00 & 2.67e-10 & 6.00 & 2.11e-10 & 5.99 \\ 
\hline 
& & & $m = 2$ & & &\\
\hline 
 10 & 3.62e-06 & - & 1.30e-06 & - & 9.47e-07 & - \\ 
20 & 1.40e-08 & 8.01 & 5.12e-09 & 7.99 & 3.88e-09 & 7.93 \\ 
40 & 5.49e-11 & 7.99 & 1.99e-11 & 8.00 & 1.51e-11 & 8.01 \\ 
80 & 2.18e-13 & 7.98 & 7.78e-14 & 8.00 & 5.97e-14 & 7.98 \\ 
\end{tabular}
} 
\end{center} 
\label{tab:Taylor2DsinMult}
\end{table} 

{}

\subsubsection{Second Half-step and Boundary Conditions}
To complete a full time step we repeat this procedure, starting with the initial data obtained from evolving \eqref{ODEs1D} a half-step $\Delta t / 2$ on the dual grid. At the interior nodes $x_i, \ i = 1, \dots, n_x-1$ the procedure is the same as above, but at the boundary nodes we must fill in ghost polynomials at $x_{i-\frac{1}{2}}$ and $x_{n_x+\frac{1}{2}}$, for example by using the properties of the PDE or as in the case of periodic boundary conditions considered here by simply copying data from the opposite boundary. 

\subsection{Hermite Method in Two-Dimensions}
The method in higher dimensions is a direct tensor product extension of the one dimensional procedure. We now describe the method for the Hamilton-Jacobi equation in two-dimensions, with periodic boundary conditions 
\begin{equation}
u_t = \begin{cases} 
\varphi_t +H(\varphi_{x},\varphi_{y}) = 0,\\
\varphi(x,y,0)=\varphi^0(x,y),\\
\end{cases} 
\label{PDE2D}
\end{equation}
and on the rectangular domain $D = [x_L,x_R] \times [y_B,y_T]$. 

\subsubsection{Initialization}
We discretize the grid as follows,
\begin{equation}
(x_i,y_j) = (x_L + ih_x,y_B+jh_y), 
\label{2Dgrid}
\end{equation}
where $h_x = (x_R - x_L)/N_x$ and $h_y = (y_T - y_B)/N_y$ where $i = 0,\dots, n_x$, $j = 0,\dots, n_y$ for the primal grid and $i = 1/2, \dots , n_x - 1/2, \ j = 1/2, \dots , n_y - 1/2$  for the dual grid.

The first step initialize the degrees of freedom by setting them to be the scaled derivatives at each point on the primal grid 
\[ c_{l_x,l_y} = \left.\frac{h_x^{l_x}}{l_x!} \frac{h_y^{l_y}}{l_y!} \frac{d^{l_x}d^{l_y}  u(x,0)}{dx^{lx}dy^{ly}} \right|_{(x=x_i,y=y_j)}. \]
Note that $c_{l_x,l_y} = c_{l_x,l_y,i,j}$, but we suppress the spatial indices for notational convenience. We interpolate onto the dual grid to obtain the tensor-product Taylor polynomials
\begin{equation}
v_{i+\frac{1}{2},j+\frac{1}{2}}(x,y,0) = \sum \limits_{l_x = 0}^{2m+1} \sum \limits_{l_y = 0}^{2m+1} d_{l_x,l_y} \left(\frac{x - x_{i+\frac{1}{2}}}{h_x }\right)^{l_x} \left(\frac{y - y_{j+\frac{1}{2}}}{h_y }\right)^{l_y},
\label{2Dinterpolant}
\end{equation}
where the polynomial interpolates the function values and the partial derivatives at the four corners of the cell. Algorithmically, forming of the interpolant is done by repeated one-dimensional interpolation. For example, we may interpolate in the $y$ direction, for the function and all the $x$ derivatives at grid points $(x_i,y_j)$ and $(x_i,y_{j+1})$ to obtain one interpolant centered at $(x_i,y_{j+1/2})$ and from $(x_{i+1},y_j)$ and $(x_{i+1},y_{j+1})$ to obtain another interpolant centered at $(x_{i+1},y_{j+1/2})$. These two polynomials of degree $m$ in $x$ and degree $2m+1$ in $y$ are then interpolated in the $x$ direction using one-dimensional interpolation. The final result is a polynomial on the form (\ref{2Dinterpolant}). 
\subsubsection{Evolution}
Similar to the one-dimensional case, we treat the coefficients in the Hermite expansions as functions of time. That is, we expand 
\begin{equation}
v_{i+1/2,j+1/2}(x,y,t) = \sum \limits_{l_x = 0}^{2m+1} \sum \limits_{l_y = 0}^{2m+1} d_{l_x,l_y}(t) \left(\frac{x - x_{i+\frac{1}{2}}}{h_x }\right)^{l_x} \left(\frac{y - y_{j+\frac{1}{2}}}{h_y }\right)^{l_y}.
\label{PDE_Hermite_expansion2D}
\end{equation}
For our PDE $\varphi_t = -H(\varphi_x,\varphi_y)$ we substitute \eqref{PDE_Hermite_expansion2D}:
\begin{align}
\frac{\partial v_{i+\frac{1}{2}}(x,y,t) }{\partial t } & =  \sum \limits_{l_x = 0}^{2m+1} \sum \limits_{l_y = 0}^{2m+1} d_{l_x,l_y}^{'}(t) \left(\frac{x - x_{i+\frac{1}{2}}}{h_x }\right)^{l_x} \left(\frac{y - y_{j+\frac{1}{2}}}{h_y }\right)^{l_y} \nonumber\\& = -H(\varphi_x,\varphi_y).
\label{PDE_Hermite_expansion_substitution2D}
\end{align}
We differentiate in \eqref{PDE_Hermite_expansion_substitution2D} $k$ times in the $x$-coordinate and $l$ times in the $y$-coordinate and evaluate at $(x,y) = (x_{i+\frac{1}{2}},y_{j+\frac{1}{2}})$ to obtain
\begin{equation}
\frac{k!}{h_x^{k}}\frac{l!}{h_y^{l}} d^{'}_{k,l} (t) = -\left. \frac{\partial^k}{\partial x^k} \frac{\partial^l}{\partial x^l} H(v_{x},v_{y}) \right |_{(x,y) = \left(x_{i + \frac{1}{2}},y_{j+\frac{1}{2}}\right)}.
\end{equation}
Similar to the one-dimensional case, the differentiation of a non-linear $H$ can spawn new terms by the product rule. We avoid this by approximating the Hamiltonian by a Taylor polynomial of degree $(2m+1) \times (2m+1)$ 
\begin{equation}
H(v_x,v_y) \approx \sum \limits_{l_x = 0}^{2m+1} \sum \limits_{l_y = 0}^{2m+1} b_{l_x,l_y}(t) \left(\frac{x - x_{i+\frac{1}{2}}}{h_x }\right)^{l_x} \left(\frac{y - y_{j+\frac{1}{2}}}{h_y }\right)^{l_y}.
\label{Taylor_expansion2D}
\end{equation}

From \eqref{PDE_Hermite_expansion_substitution2D} and \eqref{Taylor_expansion2D} we obtain the local system of ODEs
\begin{equation}
d^{'}_{l_x,l_y}(t) = b_{l_x,l_y}(t), \quad l_x = 0, \dots , 2m+1, \ l_y = 0, \dots, 2m+1.
\label{ODEs2D}
\end{equation}
We can evolve the solution of this system using, e.g., a Runge-Kutta method.

\subsubsection{Boundary Conditions for the Second Half-step}
To complete a full time step we repeat this procedure, starting with the initial data obtained from evolving \eqref{ODEs2D} a half-step on the dual grid. At the interior nodes $(x_i,y_j), \ i = 1,\dots, n_x-1, \ j = 1, \dots, n_y - 1$ the procedure is the same as above, but at the boundary nodes we must fill in ghost polynomials. In our case we use the periodic boundary conditions to fill in the ghost polynomials. 

\section{Adopting the Persson-Peraire, Klockner-Warburton-Hesthaven Sensor to Hermite Methods}
\subsection{Estimating Smoothness}
The smoothness detector used to modify our Hermite method is an adaptation of the sensor in \cite{persson2006sub}. The sensor uses orthogonal polynomials, $\lbrace \phi_n \rbrace_{n=0}^{N_p-1}$, on each cell to estimate the smoothness of the solution, where $N_p$ is the order of approximation. The method developed here also utilizes the improvements to the Persson-Peraire sensor  proposed by Klockner, Warburton and Hesthaven in \cite{klockner2011viscous}. 

On each element $D_k$, the Persson-Peraire sensor computes a smoothness indicator 
\begin{equation}
S_k = \frac{(q_N,\phi_{N_{p-1}})^2_{L^2_{(D_k)}}}{||q_N||^2_{L^2_{(D_k)}}},
\label{smoothness_sensor}
\end{equation}
where $N$ is the degree of the interpolating polynomials. \deaa{While this sensor is easy to use and implement we have found that the improved approach in \cite{klockner2011viscous} is more robust. We now explain how we use the approach of \cite{klockner2011viscous} in the context of Hermite methods.} 

The smoothness estimator we use relies on a projected version of the the solution. Precisely in one-dimension we project onto the orthogonal basis spanned by Legendre polynomials in a cell and in two dimensions we project onto a tensor product basis of one-dimensional Legendre polynomials. We take $N = 2m+1$ and \Allen{$N_p = 2m+3$}. 

For example, (in one dimension), let $q_N = \sum \limits_{n = 0}^{N_{p-1}} \hat{q}_n \phi_n$ be the projection, then its modes (coefficients) decay according to 
\begin{equation}
|\hat{q}_n| \sim cn^{-s}.
\label{modal_decay}
\end{equation}
Taking the logarithm of \eqref{modal_decay} we have 
\begin{equation}
\log |\hat{q}_n| \approx \log (c) -s \log (n),
\label{log_modal_decay}
\end{equation}
we may find the $c$ and $s$ via minimizing the least squares function 
\begin{equation}
\sum \limits_{n = 1}^{N_{p-1}} | \log |\hat{q}_n| - (\log(c) - s \log(n))|^2.
\end{equation}
Note that the sum begins at $n=1$ \Allen{(since $\log(0) = -\infty$ we minimize without $n = 0$)}, thus the constant coefficient data is not taken into account when estimating smoothness.

\Allen{The removal of the constant-mode information from the estimation process can cause problems. Consider a constant function perturbed by white noise. Since the constant-mode information is removed the smoothness detector only sees the white noise, which could lead to an erroneous smoothness estimate. A fix to this problem, a technique called \textit{baseline modal decay}, was introduced in \cite{klockner2011viscous}.} Heuristically, the idea is to re-add the sense of scale by distributing the energy according to a perfect modal \Allen{decay}
\begin{equation}
|\hat{b}_n| \sim C\frac{1}{n^N},
\end{equation} 
for $N$ the polynomial degree of the method, the normalizing factor $C$ \Allen{is chosen to enforce} 
\[
\sum \limits_{n = 1}^{N_p - 1} |\hat{b}_n|^2 = 1.
\]
We input the coefficients 
\begin{equation}
|\tilde{q}_n|^2 := |\hat{q}_n|^2 + ||q_N||_{L^2_{D_k}} |b_k|^2 \quad \text{for} \ n \in \lbrace 1, \dots , N_p - 1 \rbrace,
\end{equation}
into skyline pessimization (described below), instead of the raw coefficients $|\hat{q}_n|$.

There are certain situations where the estimator can be fooled. For example, when interpolating the function $x\Theta(x)$, where $\Theta(x)$ is the Heavyside function, \cite{klockner2011viscous} shows that there is an odd - even effect for which odd modes greater than three are numerically zero. A correction of this problem was given by introducing a technique named \textit{skyline pessimization}. The idea is as follows: if you have a mode $n$ with a small coefficient $|\hat{q}_n|$ such that there exists another coefficient with $m > n$ and $|\hat{q}_m| \gg |\hat{q}_n|$, then the coefficient $|\hat{q}_n|$ is most likely spurious and should not be taken into account when estimating $s$. Therefore, the idea is to generate a new set of modal coefficients by 
\begin{equation}
\bar{q}_n := \max \limits_{i \in \lbrace \min (n,N_p-2), \dots, N_p -1 \rbrace} |\hat{q}_i | \quad \text{for} \  n \in \lbrace{1,2,\dots,N_p - 1 \rbrace}.
\label{skyline}
\end{equation}
This forces each modal coefficient to be raised up to the largest higher-numbered modal coefficient, eliminating non-monotone decay. 

\subsection{Computing the Viscosity from the Smoothness Sensor}

\subsubsection{Implementation in One-Dimension}
\begin{figure}[hbt]
\begin{center}
\includegraphics[width=0.45\textwidth]{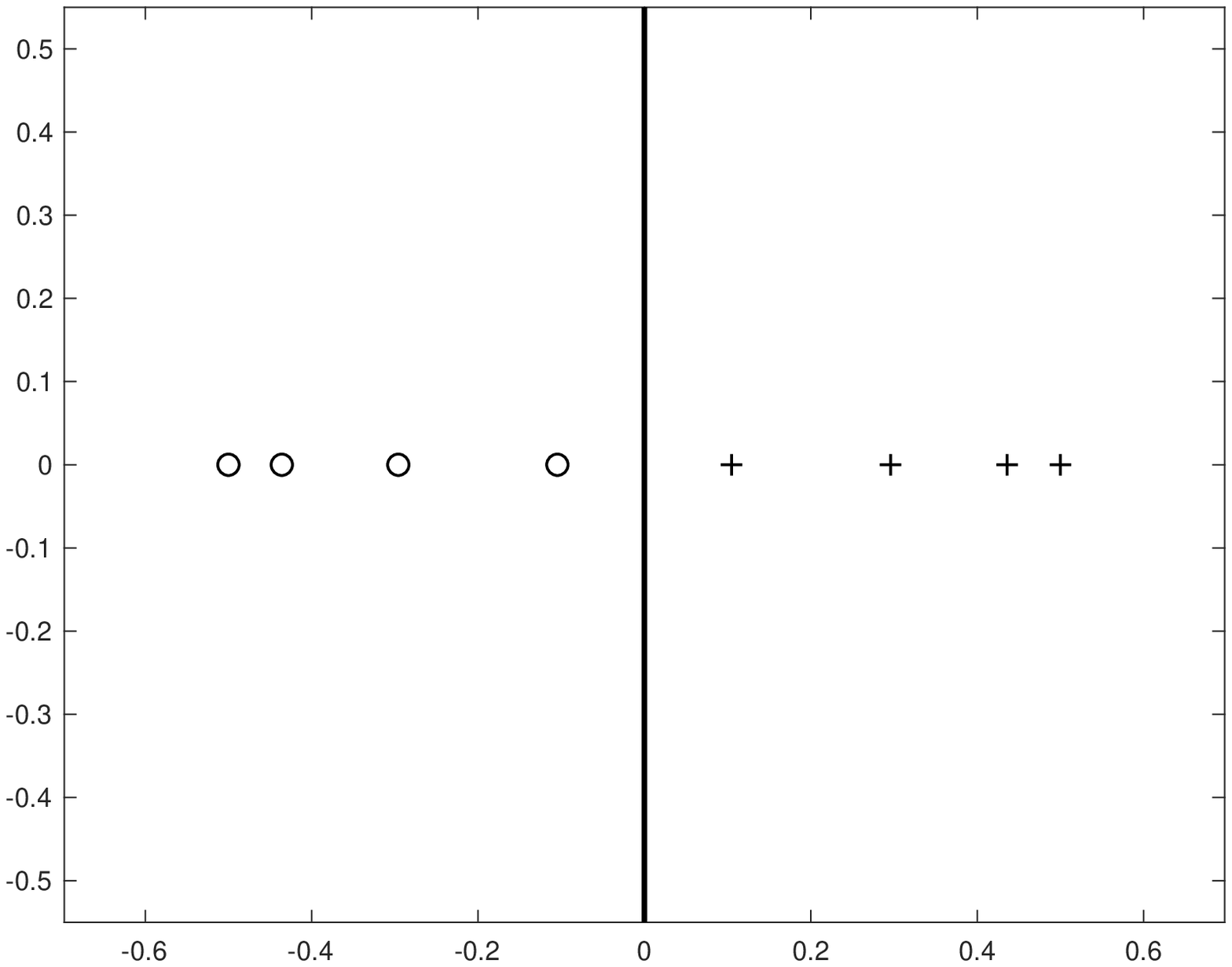}
\includegraphics[width=0.45\textwidth]{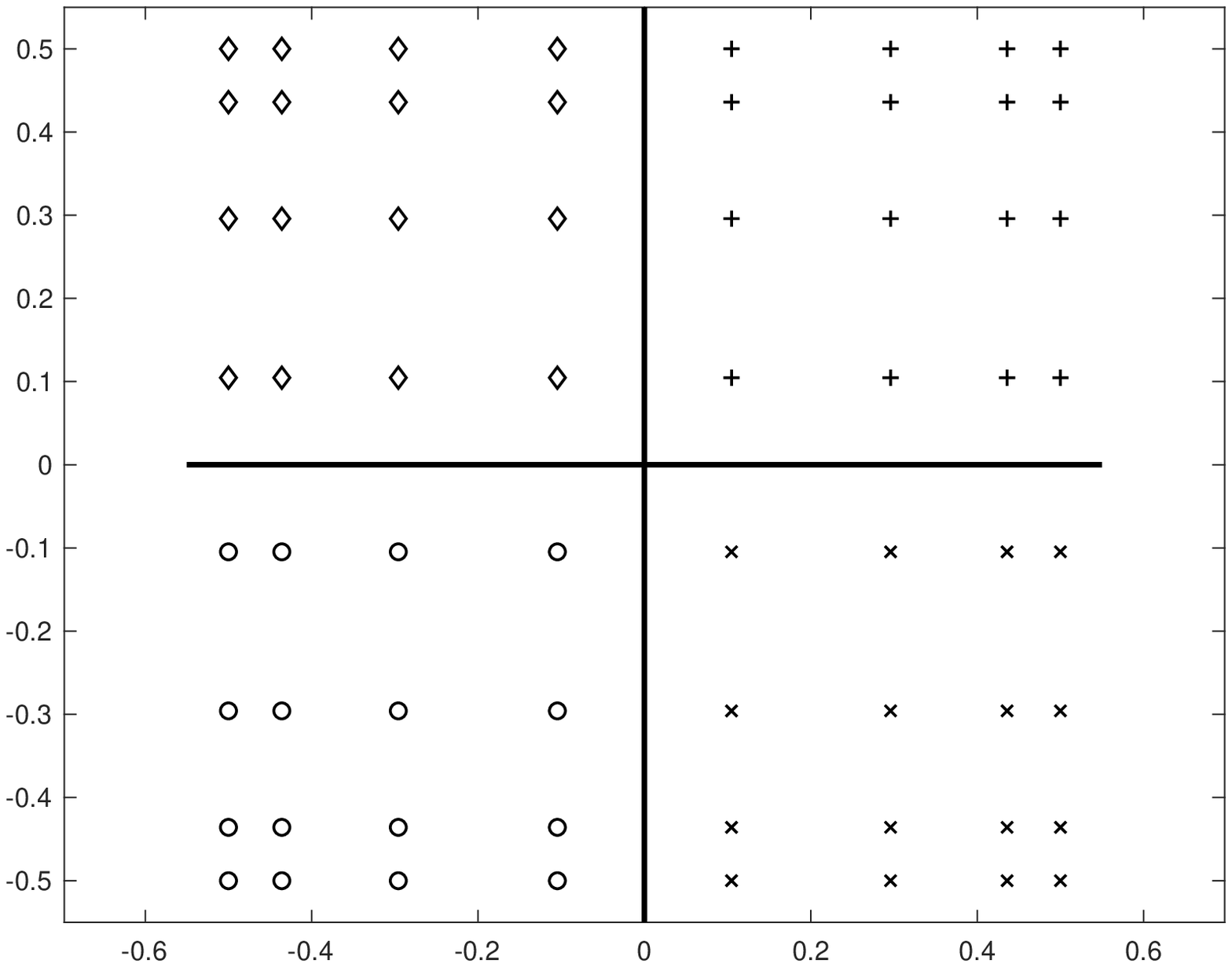} 
\caption{Left LGL nodes in one-dimension; Right LGL nodes in two-dimensions. \label{SensorFig}}
\end{center}
\end{figure}

At each time step, $t_k$, we approximate the smoothness of the function inside each cell. In order to determine the modal coefficients, $\lbrace \hat{q}_n \rbrace$, we take $2m+1$ Legendre-Gauss-Lobatto nodes (see Figure \ref{SensorFig}) $z_i$ and map the nodes in [-1,0] to $[x_k,x_{k + 1/2}]$ and the nodes in $[0,1]$ to $[x_{k + 1/2},x_{k + 1}]$, then we evaluate the Hermite interpolants $p_k(x)$ in $[x_k,x_{k + 1/2}]$ and $p_{k+1}(x)$ in $[x_{k + 1/2},x_{k + 1}]$ on the nodes to obtain the function values required for projection. Once the projection is completed, we compute $\lbrace \hat{q}_n \rbrace$ and approximate $s$ using \eqref{log_modal_decay}. \Allen{Note that the representation of the solution on each cell is a polynomial. This means that if we compute the modal coefficients using one cell, then the sensor will determine the solution is smooth (since it is sensing the smoothness of a polynomial). By taking the left and right half of each cells we are able to estimate the smoothness of the actual solution.}

If the solution is no longer $C^1$, then we introduce numerical viscosity by taking $1 - \kappa_w(s)$ where $s$ is the smoothness estimate and $\kappa_w(s)$ is given by \cite{persson2006sub} as 
\begin{equation}
\kappa_w(s) = \nu_0
\begin{cases}
0 \quad & s_w < s_0 - w,\\
\frac{1}{2} \left(1 + \sin \frac{\pi (s_w - s_0)}{2 w}\right) \quad &s_0 - w \leq s_w \leq s_0 + w,\\
1 \quad &s_0 + w < s_w,
\end{cases}
\end{equation}
where we choose $s_0 = 2$ and $s_w = 3$. These choices activate the viscosity as soon as the solution fails to be $C^1$, $s=3$, and gives maximum viscosity, $\nu_0$, when the solution is $C^0$. 

We set the maximum viscosity, $\nu_0$, to be 
\begin{equation}
\nu_0 = \lambda\frac{h}{N},
\label{nu1D}
\end{equation}
where $\lambda$ is the maximum local characteristic velocity. We estimate $\lambda$ by taking the derivative of the Hamiltonian with respect to $\varphi_x$ at each cell center and taking the maximum of the absolute value:
\begin{equation}
\lambda = \max \limits_{\Allen{i \in \lbrace 0, \dots, n_x \rbrace}} |H_{v_x}(v_x(x_{i+\frac{1}{2}}))|.
\label{lambda1D}
\end{equation}
\Allen{Here $\lambda$ is estimated using every gridpoint.} Note that this value is directly accessible, as $v_x$ is part of the degrees of freedom. 

Before modifying the PDE we average the viscosity in the spatial domain by setting 
\[\bar{\kappa}_{i} = \frac{1}{4}(\kappa_{i-1} + 2\kappa_i + \kappa_{i+1}).\]
We introduce the numerical viscosity at each timestep when evolving our PDE after interpolation. That is, after the interpolation step we modify the PDE by adding $\bar{\kappa}_i u_{xx}$ making the equation 
\begin{equation}
v_t + H(v_x) = \bar{\kappa}_i v_{xx},
\end{equation}
on the $i$th cell.
\subsubsection{Implementation in Two-Dimensions}
As our orthogonal basis we choose the tensor-product of Legendre polynomials. To check the smoothness of a cell we use the nearest interpolants on the other grid located on the four corners of the cell. That is, if the smoothness information for the cell centered at $(x_{k+1/2},y_{l+1/2})$ then we use the Hermite interpolants at $(x_k,y_l)$, $(x_{k+1},y_{l})$, $(x_{k},y_{l+1})$, and $(x_{k+1},y_{l+1})$. As in one dimension, we partition the cell into four regions and evaluate the function on each region using the interpolant that corresponds to it. For example, for the lower left region we use the information from the Hermite interpolant centered at $(x_k,y_l)$.

We analyze the decay of coefficients by grouping them by the total degree i.e., the degree of Legendre polynomial in $x$ plus the degree of the Legendre polynomial in $y$. The tensor-product Legendre polynomial is of the form
\begin{equation}
p = \sum \limits_{k = 0}^{2m+1} \sum \limits_{l = 0}^{2m+1} C_{k,l} \phi_k \phi_l,
\end{equation}
where $\phi_k$ and $\phi_l$ are the Legendre polynomials in one-dimension. 

There are several ways to order the polynomials in each degree; however, some orderings will fool the sensor. We take the maximum coefficient in absolute value for each total degree and use that as our input to the sensor for that degree. That is, we take $c_i = \max|C_{k,l}|$, where $k+l = i$. Once we obtain this ordering we apply \textit{baseline modal decay} and \textit{skyline pessimization} in the same way as the one-dimensional case.

We set the maximum viscosity, $\nu_0$, to be 
\begin{equation}
\nu_0 = \lambda\frac{h}{N},
\label{nu2D}
\end{equation}
where $N$ is the degree of the Hermite interpolating polynomial in one coordinate direction and $h = \max\lbrace h_x,h_y\rbrace$. For estimating $\lambda$ we adapt the Lax-Friedrichs flux given in \cite{crandall1980monotone,osher1991high} by taking the partial derivatives of the Hamiltonian with respect to $\varphi_x$ and $\varphi_y$ evaluating them at the cell centers and taking the maximum of the absolute value:
\begin{equation}
\lambda = \max \limits_{i,j}\lbrace |H_{v_x}(v_x(x_{i+\frac{1}{2}},y_{j+\frac{1}{2}}),v_y(x_{i+\frac{1}{2}},y_{j+\frac{1}{2}}))|,|H_{v_y}(v_x(x_{i+\frac{1}{2}},y_{j+\frac{1}{2}}),v_y(x_{i+\frac{1}{2}},y_{j+\frac{1}{2}}))| \rbrace.
\label{lambda2D}
\end{equation}
Before modifying the PDE we average the viscosity in the spatial domain by setting 
\begin{align*}\bar{\kappa}_{k,l} & = \frac{1}{16}(\kappa_{k-1,l-1}+\kappa_{k+1,l-1}+\kappa_{k-1,l+1}+\kappa_{k+1,l+1}+\\
 &2(\kappa_{k,l-1}+\kappa_{k-1,l}+\kappa_{k+1,l}+\kappa_{k,l+1})+4\kappa_{k,l}).
\end{align*}
We introduce the numerical viscosity at each timestep when evolving our PDE after interpolation. That is, after the interpolation step we modify the PDE making the equation 
\begin{equation}
v_t + H(v_x,v_y) = \bar{\kappa}_{k,l}(v_{xx} + v_{yy}).
\end{equation}

\begin{figure}[hbt]
\begin{center}
\includegraphics[width=0.24\textwidth]{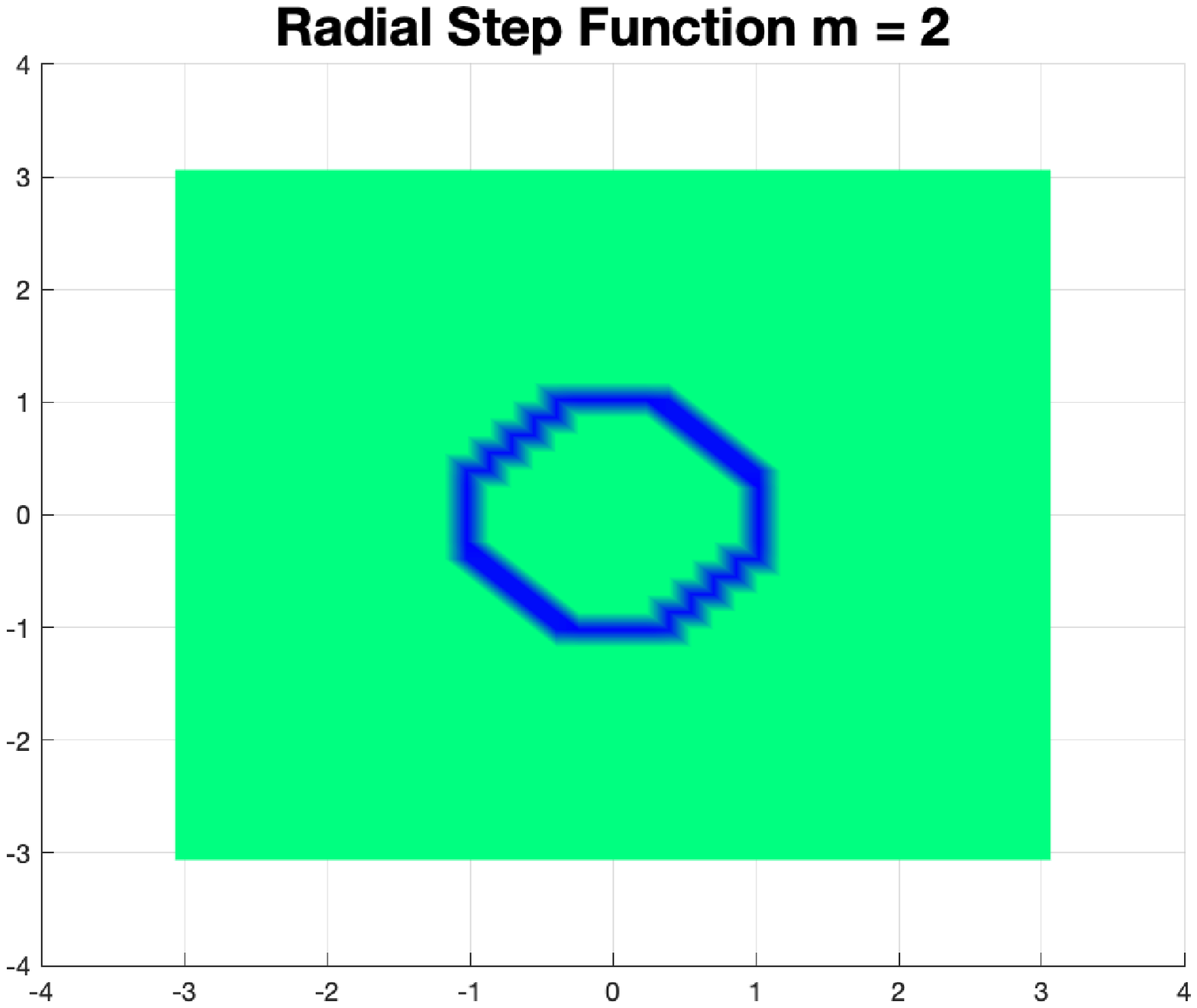} 
\includegraphics[width=0.24\textwidth]{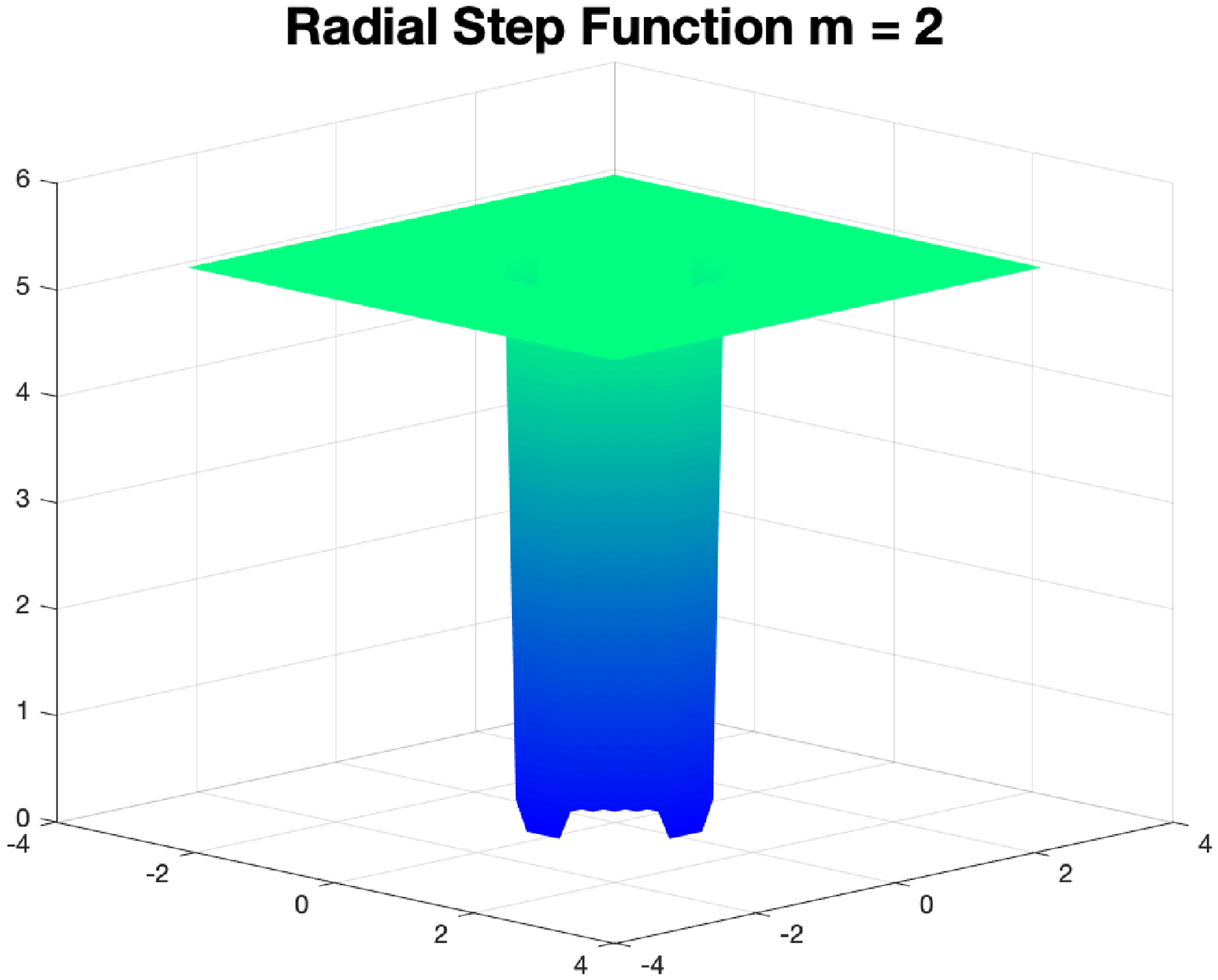} 
\includegraphics[width=0.24\textwidth]{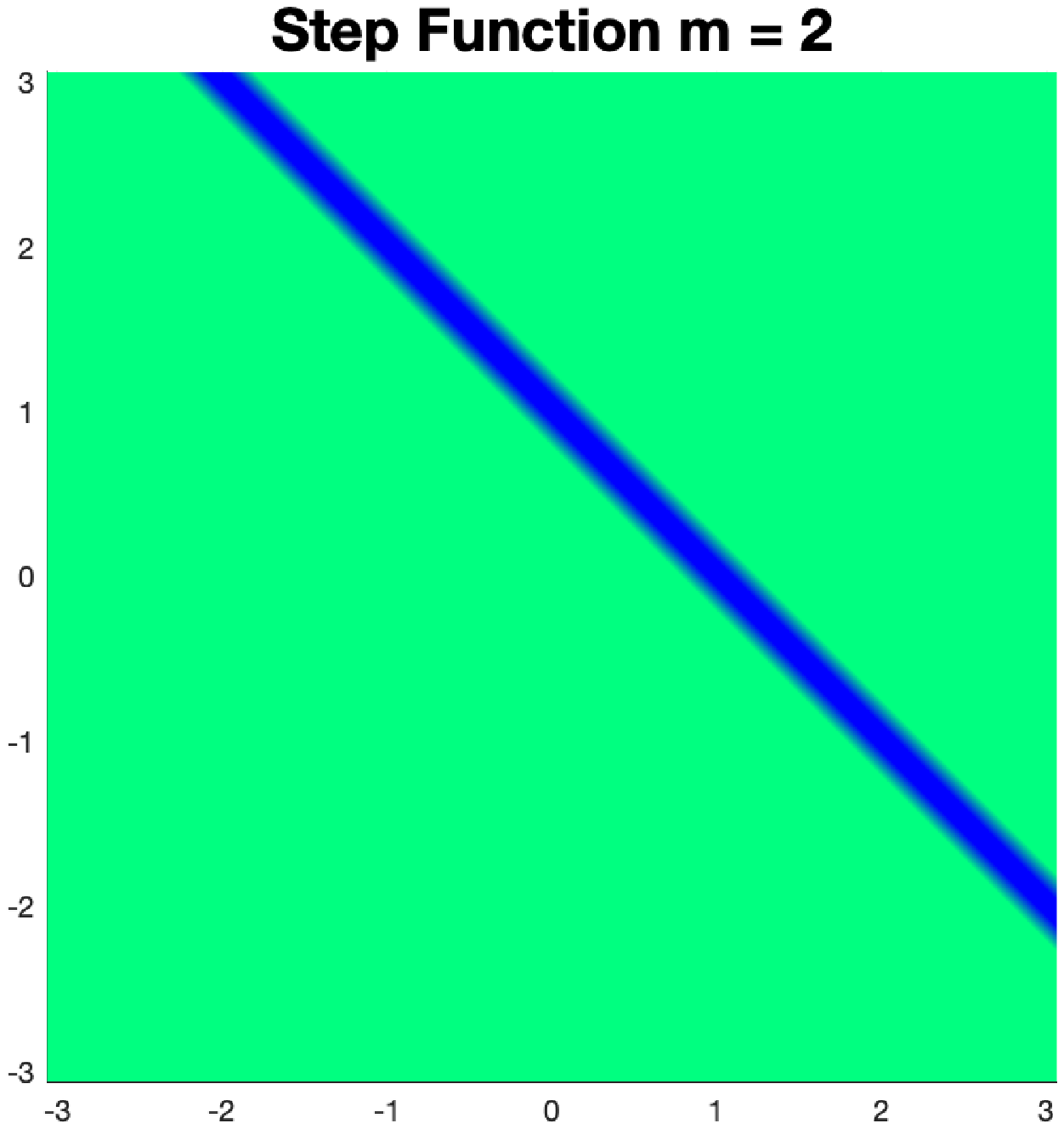} 
\includegraphics[width=0.24\textwidth]{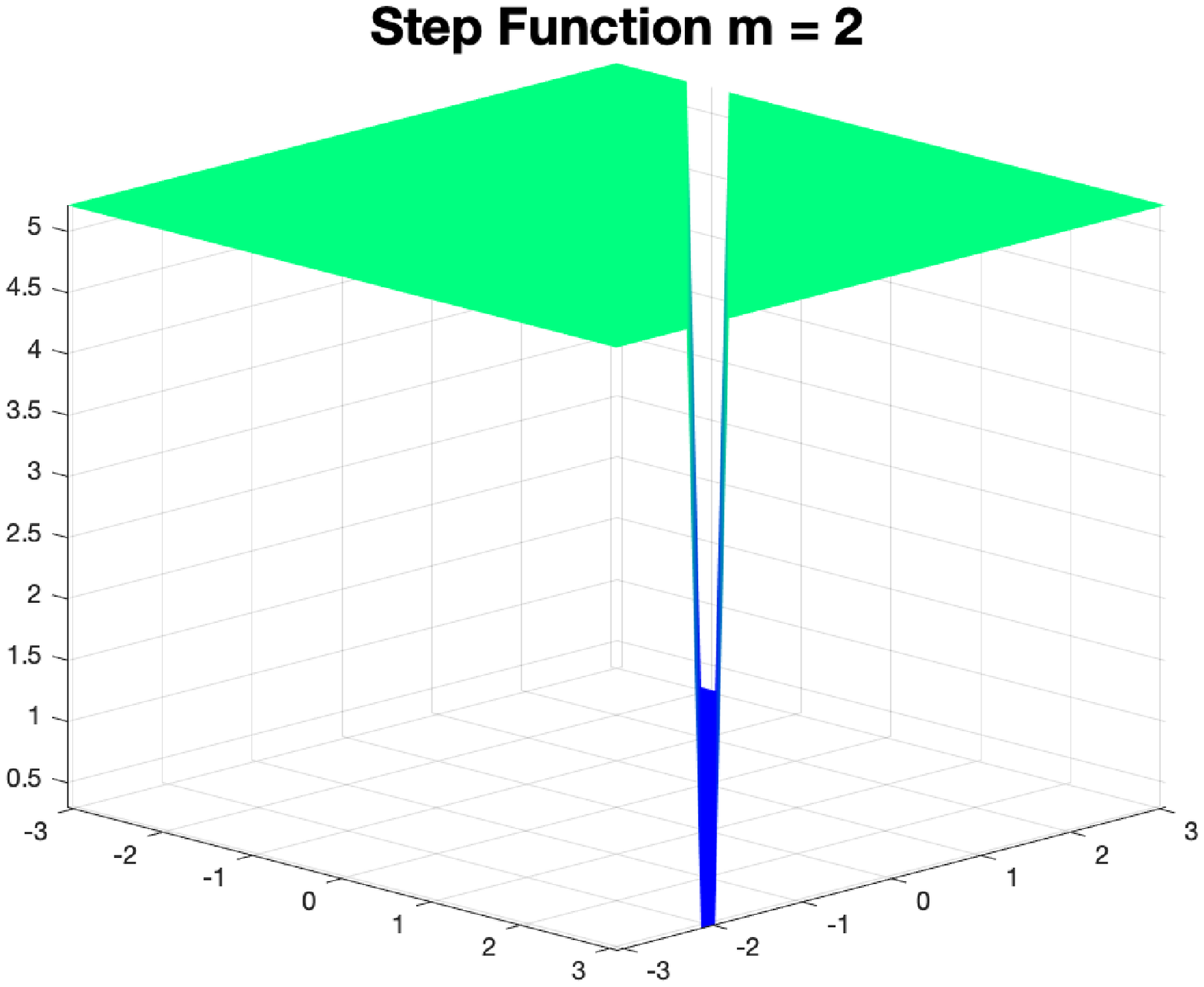} 
\caption{Top: different views of the smoothness sensor when applied to a radial stepfunction $f(x) = 1$ if $r^2 \leq 1$ and zero otherwise. In the smooth region the smoothness is estimated to be 5 and near the the discontinuity it is estimated to be approximately 1/2. Bottom: different views of the smoothness sensor for a stepfunction with $f(x) = 1$ if $x+y\leq 1$ and zero otherwise. The results are similar to the radial case. \label{fig:smoothess_test}}
\end{center}
\end{figure}

In Figure \ref{fig:smoothess_test} we test the smoothness sensor on two discontinuous functions, a radially symmetric step function and an oblique step function. In both cases the sensor correctly  determines the level of smoothness of the underlying function. 

\section{Numerical Examples}
For each of the following examples we expect the convergence rate to be $2m+1$ when the solution is smooth, where $m$ is number of derivatives used in the interpolation process. When the solution fails to remain smooth we do not expect to see the optimal convergence rates, instead we seek sharp resolution of kinks. For each example we give convergence rates at a time when the solution is still smooth. For the timestepping we use the classical Runge-Kutta method of order 4 (RK4). To this end we chose the timestep small enough so that the temporal error is dominated by the spatial error. \FinalEdits{The derivation of analytical solutions for examples where convergence rates are estimated can be found in the appendix.}
 
\subsection{Examples in One Dimension} 
\subsubsection{Example 1}
\begin{figure}[htb]
\begin{center}
\includegraphics[width=0.6\textwidth]{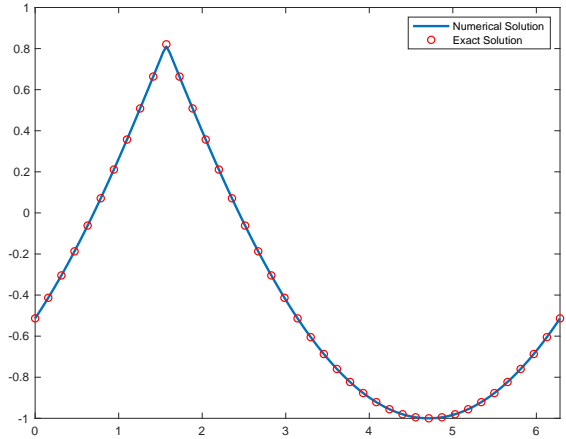} 
\caption{Example 1 at time $t = 1.5$ well after the kink has developed. This computation was done with $m = 3$ and $N = 80$ cells. The solid line is the numerical solution and the dots are the exact solution.\label{fig:Burgers1Dplot}}
\end{center}
\end{figure}

In this example we solve the one-dimensional Burgers' equation 
\begin{equation} 
\varphi_t +\frac{1}{2} \Allen{(\varphi_x)^2}=0,
\nonumber
\end{equation}
with initial condition \mbox{$\varphi(x,0)=\sin(x)$} and with periodic boundary conditions \mbox{$\varphi(0,t)=\varphi(2\pi,t)$}. The solution is smooth until time $T = 1.0$, at this time a shock will form in $\varphi_x$. Our grid is determined by $x_l = 0$, $x_r = 2\pi$ and the number of grid points, $n_x$. For this example we start with $n_x = 20$ and refine the grid by a factor of two until $n_x = 160$ in order to demonstrate convergence to the viscosity solution. Before the solution develops a kink we demonstrate that our method achieves $2m+1$ order accuracy at time $T=0.5$ as evidenced by the errors measured in the  $L_1, L_2$ and $L_{\infty}$ norm along with the estimated rates of convergence reported in Table \ref{tab:Example1smooth}. We also demonstrate that we converge to the viscosity solution at time $T = 1.5$ in Figure \ref{fig:Burgers1Dplot} along with the errors  reported in Table \ref{tab:Example1Notsmooth}. \FinalEdits{Note that the kink formed closely resembles an absolute value function. The degree $2m+1$ Hermite interpolant of the absolute value function, $|x|$, can be explicitly written down. It is  
\begin{equation*}
p(x) = \sum \limits_{k = 0}^{m} \binom{2k}{k} \frac{(-1)^{k+1}(x^2-1)^k}{2^{2k}(2k-1)},
\end{equation*}
and this in turn corresponds to the first terms of the generalized binomial expansion  
\[
(1+t)^{\frac{1}{2}} = \sum \limits_{k = 0}^{\infty} \binom{1/2}{k}t^k,
\]
which is a convergent approximation of the absolute value function when replacing $t = x^2 - 1$ and when $|x|<1$. At $x=0$ and for a fixed $m$ this approximation is positive and there is thus a $\mathcal{O}(1)$ error in a single cell of width $h$. This is the source of the $\mathcal{O}(h)$ in the max norm and $\mathcal{O}(h^2)$ in the 2 norm.}

\begin{table}[ht]
 \caption{Errors in Example 1 in the $L_1$, $L_2$ and $L_{\infty}$ norms at time $T = 1.5$ are displayed along with estimated rates of convergence. Note that these errors occur after a kink has formed.} 
\begin{center} 
 \begin{tabular}{c c c c c c c} 
 \hline
 n & $L_1$ error & Conv. Rate & $L_2$ error & Conv. Rate & $L_{\infty}$ error & Conv. Rate \\ 
\hline 
& & & $m = 2$ & & &\\
 \hline 
 20 & 2.76e-05 & - & 2.70e-05 & - & 5.47e-05 & - \\ 
40 & 9.45e-07 & 4.87 & 9.53e-07 & 4.83 & 2.03e-06 & 4.75 \\ 
80 & 3.05e-08 & 4.95 & 3.10e-08 & 4.94 & 6.59e-08 & 4.95 \\ 
160 & 9.42e-10 & 5.02 & 9.56e-10 & 5.02 & 2.01e-09 & 5.03 \\ 
\hline 
& & & $m = 3$ & & &\\
 \hline 
 20 & 3.56e-07 & - & 4.55e-07 & - & 1.26e-06 & - \\ 
40 & 2.61e-09 & 7.09 & 3.22e-09 & 7.14 & 8.49e-09 & 7.21 \\ 
80 & 2.03e-11 & 7.01 & 2.50e-11 & 7.01 & 6.23e-11 & 7.09 \\ 
160 & 1.82e-13 & 6.80 & 1.85e-13 & 7.08 & 4.56e-13 & 7.10 \\ 
\hline 
\end{tabular} 
\end{center} 
\label{tab:Example1smooth}
\end{table} 

\begin{table}[ht]
 \caption{Errors in Example 1 in the $L_1$, $L_2$ and $L_{\infty}$ norms at time $T = 1.5$ are displayed along with estimated rates of convergence. Note that these errors occur after a kink has formed.} 
\begin{center} 
 \begin{tabular}{c c c c c c c} 
 \hline
 n & $L_1$ error & Conv. Rate & $L_2$ error & Conv. Rate & $L_{\infty}$ error & Conv. Rate \\ 
\hline 
& & & $m = 2$ & & &\\
 \hline 
 20 & 1.16e-02 & - & 1.71e-02 & - & 4.00e-02 & - \\ 
40 & 2.73e-03 & 2.09 & 5.86e-03 & 1.54 & 1.97e-02 & 1.02 \\ 
80 & 6.81e-04 & 2.00 & 2.08e-03 & 1.50 & 9.85e-03 & 1.00 \\ 
160 & 1.70e-04 & 2.00 & 7.28e-04 & 1.51 & 4.87e-03 & 1.01 \\ 
\hline 
& & & $m = 3$ & & &\\
 \hline 
  20 & 6.87e-03 & - & 1.46e-02 & - & 3.67e-02 & - \\ 
40 & 1.66e-03 & 2.05 & 4.94e-03 & 1.57 & 1.75e-02 & 1.07 \\ 
80 & 4.12e-04 & 2.01 & 1.75e-03 & 1.50 & 8.75e-03 & 1.00 \\ 
160 & 1.02e-04 & 2.01 & 6.18e-04 & 1.50 & 4.38e-03 & 1.00 \\ 
\hline 
\end{tabular} 
\end{center} 
\label{tab:Example1Notsmooth}
\end{table} 

\begin{figure}[htb]
\begin{center}
\includegraphics[width=0.48\textwidth]{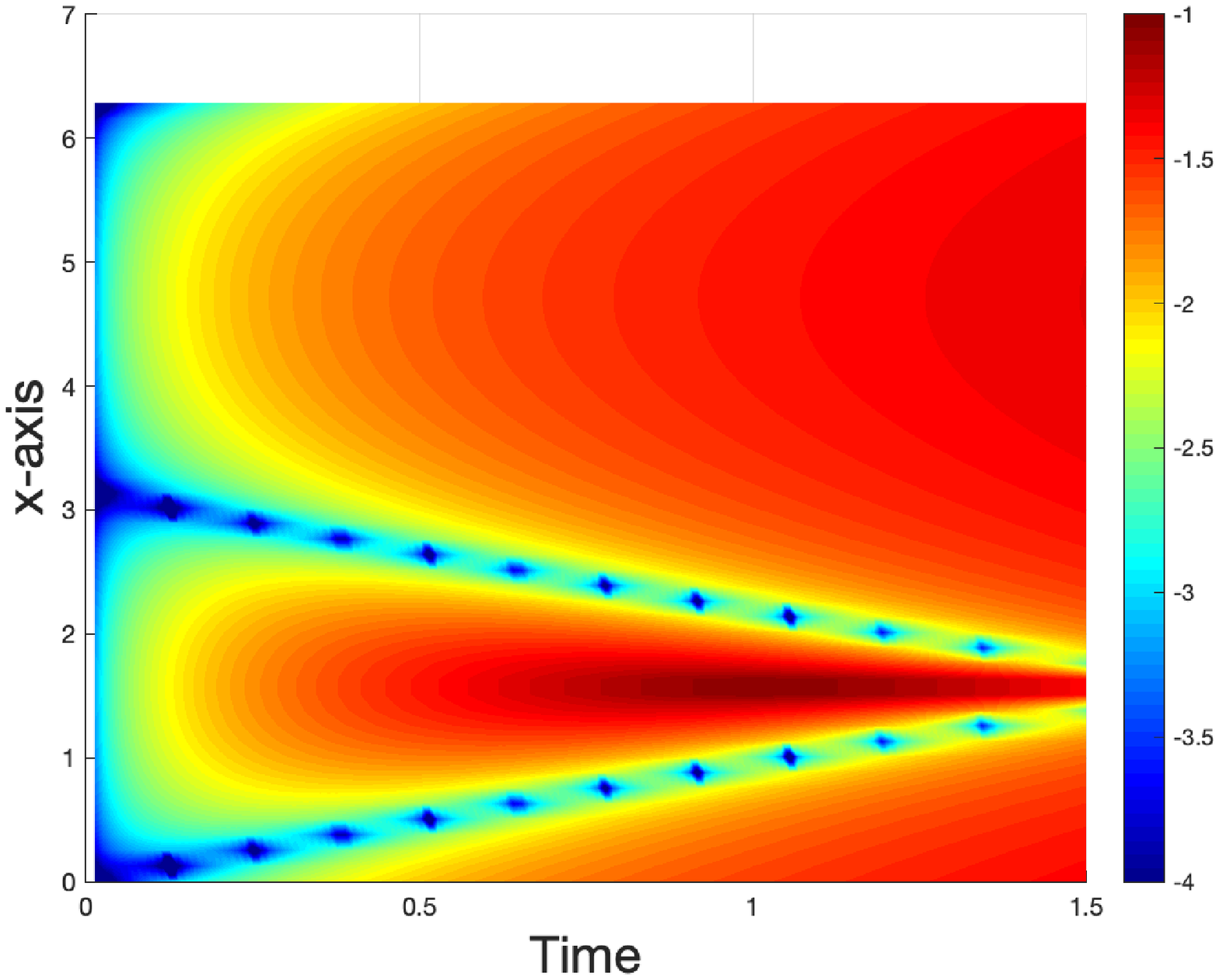}
\includegraphics[width=0.48\textwidth]{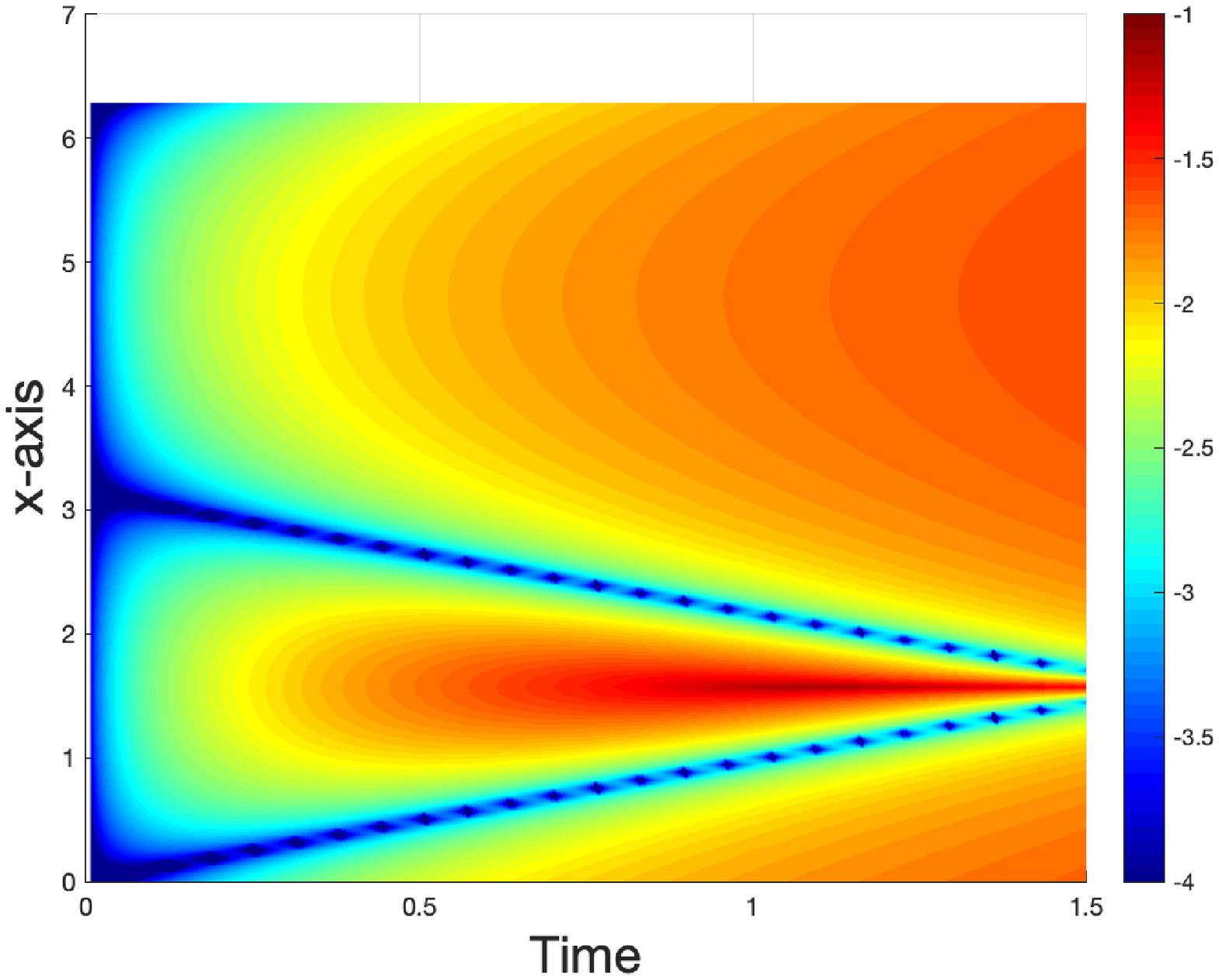}
\caption{Here we display the evolution of the errors with a refinement. On the left we display the evolution of errors with 50 cells and on the right we display the evolution of errors with 100 cells. \label{fig:spaceTimeBurgers}}.
\end{center}
\end{figure}

The convergence rates displayed in the tables show us that we are converging to the viscosity solution as the grid is being refined. We observe in Figure \ref{fig:Burgers1Dplot} that the method is able to capture the kink formed at $\frac{\pi}{2}$. In Figure \ref{fig:spaceTimeBurgers} we see that as we refine the grid the error is localized where the kink is formed.   

\begin{figure}[ht]
\begin{center}
\includegraphics[width=0.6\textwidth]{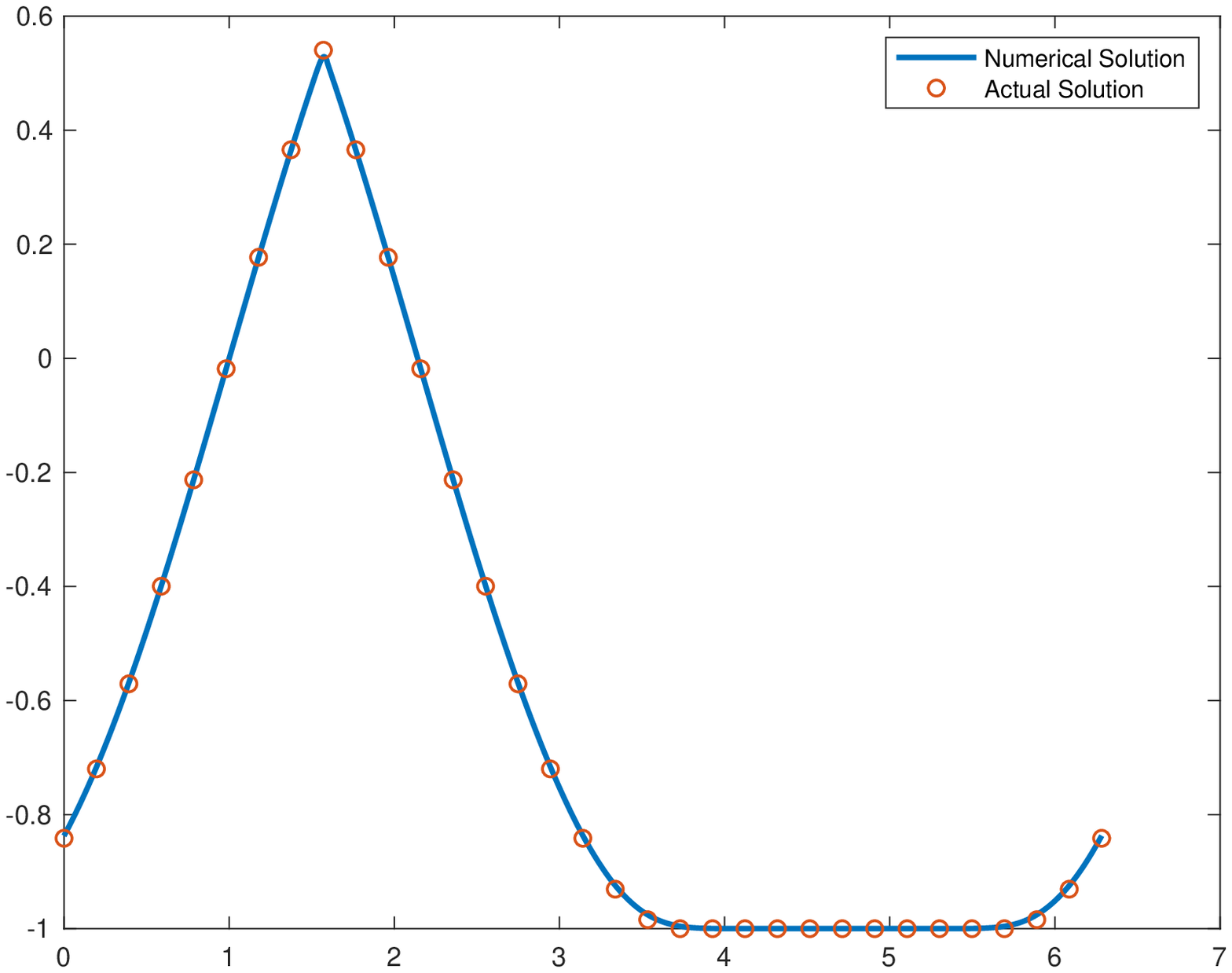} 
\caption{Example 2 at time $t = 1.0$ with $m = 3$ and $N = 80$ cells the solid line is the numerical solution and the dots are exact solution. \label{fig:Eikonal1Dplot}}
\end{center}
\end{figure}
\subsubsection{Example 2}
In this example we solve the one-dimensional Eikonal equation
\begin{align*} 
\varphi_t + |\varphi_x|=0,
\end{align*}
with initial condition \mbox{$\varphi(x,0)=\sin(x)$} and with periodic boundary conditions. The viscosity solution to this equation has a shock forming in $\varphi_x$ at $x = \pi/2$ and a rarefaction wave at $x = 3\pi/2$. The solution is nonsmooth for all $T >0$ so we do not expect order $2m+1$ convergence. To analyze the convergence we use the same grids as in Example 1. We report the $L_1,L_2$ and $L_{\infty}$ errors and their estimated rates of convergence in Table \ref{tab:Eikonal}. In Figure \ref{fig:Eikonal1Dplot} we plot the numerical solution to demonstrate convergence to the viscosity solution.
\begin{table}[ht]
 \caption{Errors in Example 2 in the $L_1$, $L_2$ and $L_{\infty}$ norms at time $T = 1.0$ are displayed along with estimated rates of convergence.} 
\begin{center} 
 \begin{tabular}{c c c c c c c} 
\hline 
  n & $L_1$ error & Conv. Rate & $L_2$ error & Conv. Rate & $L_{\infty}$ error & Conv. Rate \\ 
\hline 
& & & $m = 2$ & & &\\
\hline
 20 & 4.84e-01 & - & 2.20e-01 & - & 1.94e-01 & - \\ 
40 & 2.35e-01 & 1.04 & 1.13e-01 & 0.96 & 1.09e-01 & 0.83 \\ 
80 & 1.14e-01 & 1.04 & 5.77e-02 & 0.97 & 5.79e-02 & 0.91 \\ 
160 & 5.60e-02 & 1.03 & 2.93e-02 & 0.98 & 2.98e-02 & 0.96 \\ 
\hline 
& & & $m = 3$ & & &\\
\hline
 20 & 3.40e-01 & - & 1.58e-01 & - & 1.39e-01 & - \\ 
40 & 1.65e-01 & 1.04 & 8.11e-02 & 0.97 & 7.61e-02 & 0.87 \\ 
80 & 8.05e-02 & 1.04 & 4.13e-02 & 0.97 & 3.96e-02 & 0.94 \\ 
160 & 3.97e-02 & 1.02 & 2.10e-02 & 0.98 & 2.03e-02 & 0.97 \\ 
\hline 
\end{tabular} 
\end{center} 
\label{tab:Eikonal}
\end{table} 

\begin{figure}[ht]
\begin{center}
\includegraphics[width=0.48\textwidth]{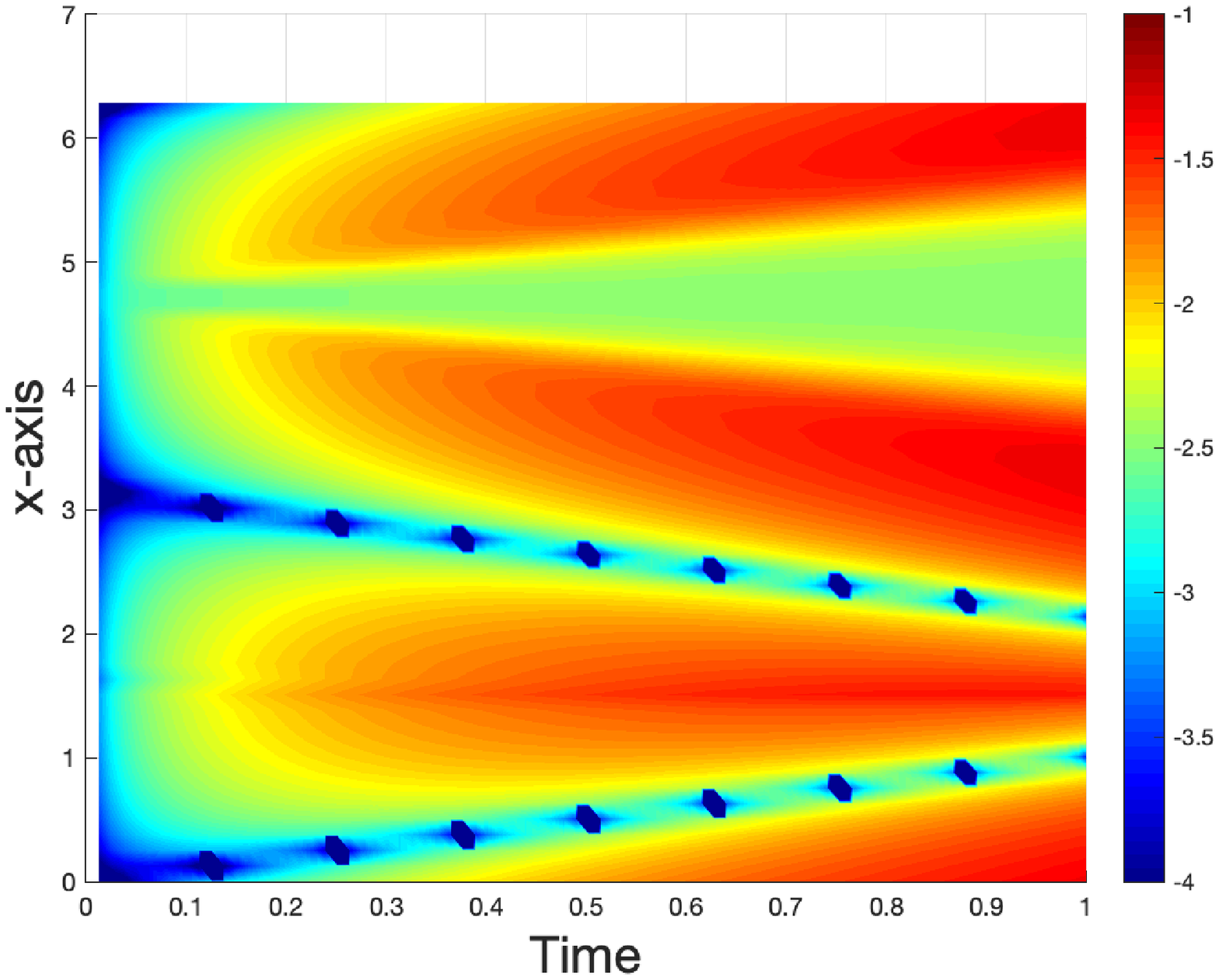}
\includegraphics[width=0.48\textwidth]{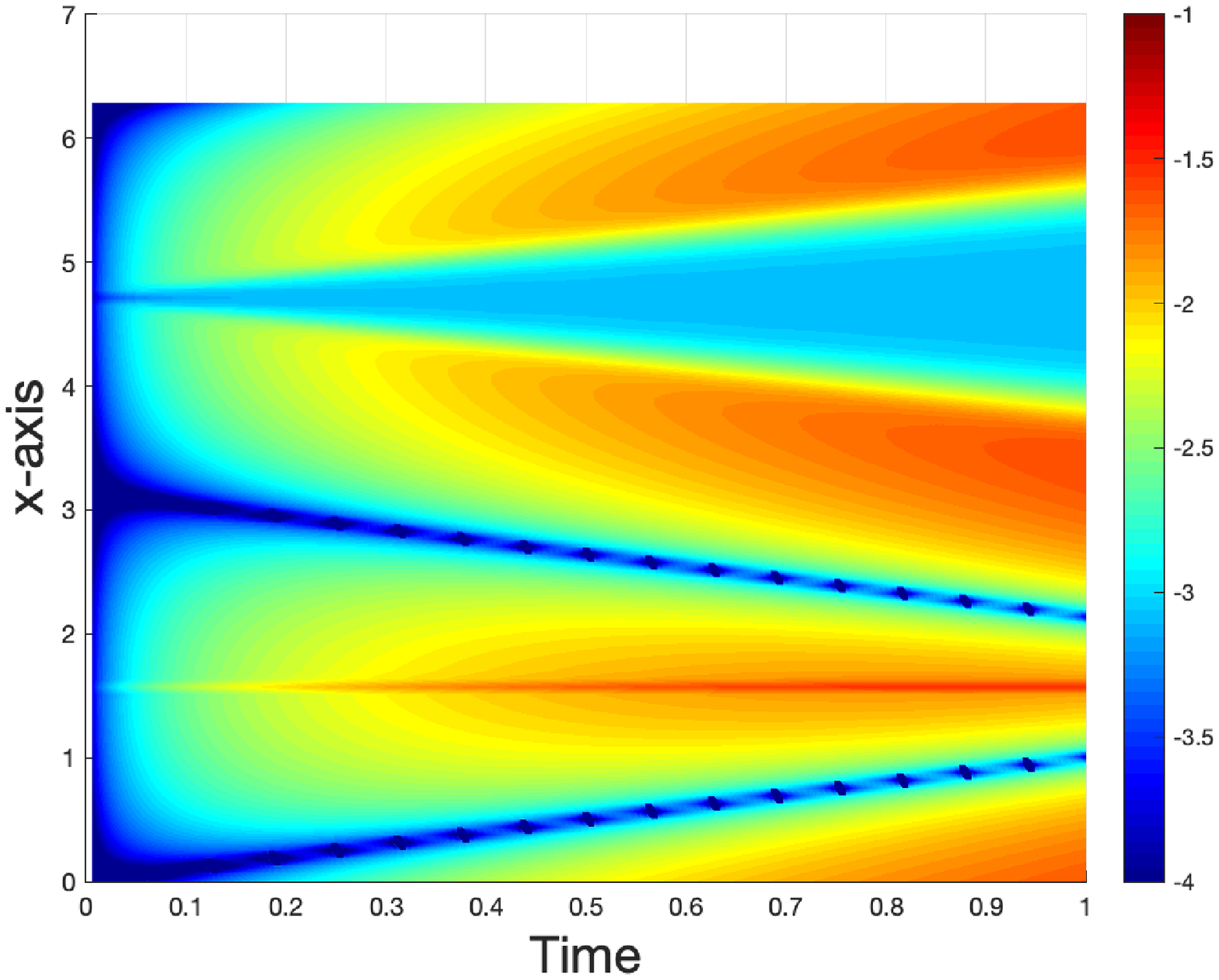}
\caption{Here we display the evolution of the errors with a refinement. On the left display the evolution of errors with 50 cells and on the right we display the evolution of errors with 100 cells. \label{fig:spaceTimeEikonal}}
\end{center}
\end{figure}

The convergence rates displayed in the tables show us that we are converging to the viscosity solution as the grid is being refined. We observe in Figure \ref{fig:Eikonal1Dplot} that the method is able to capture the kink formed at $\frac{\pi}{2}$ and the rarefaction wave at $\frac{3\pi}{2}$. In Figure \ref{fig:spaceTimeEikonal} we see that as we refine the grid the error is localized where the kink is formed. We briefly note that the discontinuity in the Hamiltonian, $H$, causes the piece-wise interpolant to lose smoothness in between cells where the sign of $\varphi_x$ changes. We plan to see if we can rectify this by using a flux-conservative formulation of this method.  

\subsubsection{Example 3}
\begin{table}[ht]
 \caption{Errors in Example 3 in the $L_1$, $L_2$ and $L_{\infty}$ norms at time $T = 0.5/{ \Allen{\pi}}^2$ are displayed along with estimated rates of convergence.} 
\begin{center} 
 \begin{tabular}{c c c c c c c } 
\hline 
& & & $m = 2$ & & &\\
\hline
  n & $L_1$ error & Conv. Rate & $L_2$ error & Conv. Rate & $L_{\infty}$ error & Conv. Rate \\  
 \hline 
 20 & 1.72e-05 & - & 3.49e-05 & - & 1.59e-04 & - \\ 
40 & 4.67e-07 & 5.20 & 1.05e-06 & 5.06 & 6.47e-06 & 4.62 \\ 
80 & 1.48e-08 & 4.98 & 2.63e-08 & 5.32 & 1.68e-07 & 5.27 \\ 
160 & 6.56e-10 & 4.50 & 8.62e-10 & 4.93 & 3.79e-09 & 5.47 \\ 
\hline 
& & & $m = 3$ & & &\\
 \hline 
 20 & 4.96e-06 & - & 1.18e-05 & - & 5.77e-05 & - \\ 
40 & 4.00e-08 & 6.95 & 1.31e-07 & 6.50 & 8.79e-07 & 6.04 \\ 
80 & 2.23e-10 & 7.48 & 7.28e-10 & 7.49 & 4.75e-09 & 7.53 \\ 
160 & 1.94e-12 & 6.84 & 4.30e-12 & 7.40 & 3.17e-11 & 7.23 \\ 
\hline 
\end{tabular} 
\end{center} 
\label{tab:nonConvexSmooth}
\end{table} 

\begin{figure}[ht]
\begin{center}
\includegraphics[width=0.6\textwidth]{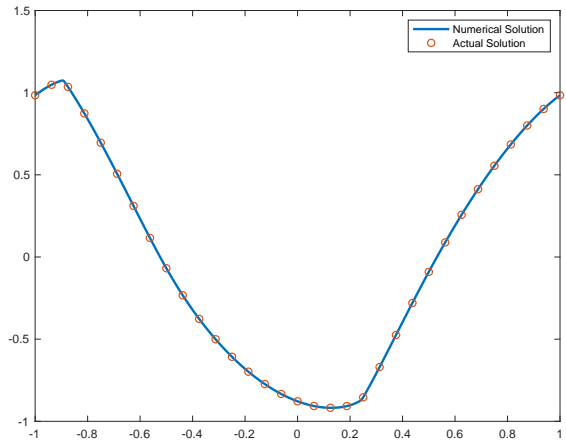}  
\caption{Example 3 at time $t = 1.5/\pi^2$ with $m = 3$ using $N = 80$ cells. The solid line is the numerical solution and the dots are the exact solution \label{fig:NonConvex1Dplot}}
\end{center}
\end{figure}

In this example we solve a one-dimensional equation with a nonconvex Hamiltonian
\begin{align*} 
\varphi_t - \cos(\varphi_x + 1) = 0,
\end{align*}
with initial condition $\varphi(x,0)=-\cos(\pi x)$ and periodic boundary conditions $\varphi(-1,t)=\varphi(1,t)$. This example shows that our scheme has high-order accuracy even when the Hamiltonian is not convex. Our grid is determined by $x_l = -1$, $x_r = 1$ and the number of grid points, $n_x$. For this example we start with $n_x = 20$ and refine the grid until $n_x = 160$ in order to demonstrate convergence to the viscosity solution. Before the solution develops a kink we demonstrate our method achieves $2m+1$ order accuracy at time $T=\frac{0.5}{\pi^2}$ by giving the $L_1, L_2$ and $L_{\infty}$ norm along with the estimated rates of convergence in Table \ref{tab:nonConvexSmooth}. We also demonstrate that we converge to the viscosity solution at time $T = \frac{1.5}{\pi^2}$ in Figure \ref{fig:NonConvex1Dplot}.

The convergence rates displayed in the tables show us that we are converging to the viscosity solution as the grid is being refined. We observe in Figure \ref{fig:NonConvex1Dplot} that the method is able to capture the kinks formed in this example. 

\Allen{
\subsubsection{Example 4}
In this example we solve a one-dimensional Riemann problem with a nonconvex Hamiltonian
\begin{align*} 
\varphi_t + \frac{1}{4}(\varphi_x^2 - 1)(\varphi_x^2 - 4) = 0,
\end{align*}
with initial data $\varphi(x,0) = -2|x|$. In this example there are two shocks propagating to the left and right connected in between by a rarefaction wave. Our grid is determined by $x_l = -1$, $x_r = 1$ and the number of grid points, $n_x$. For this example we first start with an odd number of grid points $n_x = 21$ and refine by a factor of two until $n_x = 321$ in order to demonstrate convergence to the viscosity solution. We also refine the grid for an even number of grid points starting with $n_x = 20$ and using $n_x = 320$ on the finest grid. The difference between the two is when using an even number of gridpoints the discontinuity in the initial data at $x = 0$ is located on a gridpoint $x_{n/2}$. When the discontinuity in the initial data is located on a gridpoint, $x_*$, we define the degrees of freedom at $x_*$ in two ways, depending on whether the data is being used to interpolate the left or right of $x_*$. That is, when interpolating to the dual node located to the left of $x_*$ we compute the degrees of freedom at $x_*$ using the limit of the initial data from the left, $\lim \limits_{x \to x_{*}^{-}} \varphi$, and when interpolating to the dual node located to the right of $x_*$ we compute the degrees of freedom at $x_*$ using the limit of the initial data from the right, $\lim \limits_{x \to x_{*}^{+}} \varphi$. Exact solutions are used as the boundary condition. We report the $L_1$, $L_2$ and $L_{\infty}$ norm errors at time $T = 1.0$ along with estimated rates of convergence in Table \ref{tab:Riemann1D}. We observe first order convergence for both the even and odd refinements. We note that for this more difficult example the amount of viscosity needed appears to scale inversely with the square of the CFL number. That is, at large CFL numbers we need less viscosity to compute the correct solution.

\begin{table}[ht]
 \caption{\Allen{Errors in Example 4 in the $L_1$, $L_2$ and $L_{\infty}$ norms at time $T = 1.0$ are displayed along with estimated rates of convergence.}} 
\begin{center} 
\Allen{
 \begin{tabular}{c c c c c c c } 
\hline 
 n & $L_1$ error & Convergence & $L_2$ error & Convergence & $L_{\infty}$ error & Convergence \\ 
 \hline 
 & & & Odd & & &\\
 \hline
& & & $m = 2$ & & &\\
\hline
  41 & 4.03e-02 & - & 3.83e-02 & - & 4.35e-02 & - \\
81 & 1.99e-02 & 1.02 & 1.91e-02 & 1.00 & 2.15e-02 & 1.02 \\
161 & 9.87e-03 & 1.01 & 9.54e-03 & 1.00 & 1.04e-02 & 1.05 \\
321 & 4.76e-03 & 1.05 & 4.66e-03 & 1.03 & 4.84e-03 & 1.10 \\
\hline 
& & & $m = 3$ & & &\\
\hline 
 41 & 4.92e-02 & - & 8.23e-02 & - & 6.44e-01 & - \\ 
81 & 1.97e-02 & 1.32 & 1.89e-02 & 2.12 & 2.17e-02 & 4.89 \\ 
161 & 9.79e-03 & 1.01 & 9.47e-03 & 1.00 & 1.10e-02 & 0.99 \\ 
321 & 4.87e-03 & 1.01 & 4.73e-03 & 1.00 & 5.30e-03 & 1.05 \\ 
\hline 
& & & Even & & &\\
\hline 
& & & $m = 2$ & & &\\
\hline 
40 & 3.98e-02 & - & 3.80e-02 & - & 4.32e-02 & - \\
80 & 1.95e-02 & 1.03 & 1.87e-02 & 1.02 & 2.15e-02 & 1.01 \\
160 & 9.59e-03 & 1.02 & 9.29e-03 & 1.01 & 1.02e-02 & 1.07 \\
320 & 4.62e-03 & 1.05 & 4.52e-03 & 1.04 & 4.71e-03 & 1.12 \\
\hline
& & & $m = 3$ & & &\\
\hline 
40 & 3.95e-02 & - & 3.77e-02 & - & 4.43e-02 & - \\
80 & 1.93e-02 & 1.04 & 1.86e-02 & 1.02 & 2.21e-02 & 1.01 \\
160 & 9.50e-03 & 1.02 & 9.21e-03 & 1.01 & 1.07e-02 & 1.04 \\
320 & 4.73e-03 & 1.01 & 4.59e-03 & 1.01 & 5.77e-03 & 0.89 \\
\end{tabular} 
}
\end{center} 
\label{tab:Riemann1D}
\end{table} 

\begin{figure}[ht]
\begin{center}
\includegraphics[width=0.6\textwidth]{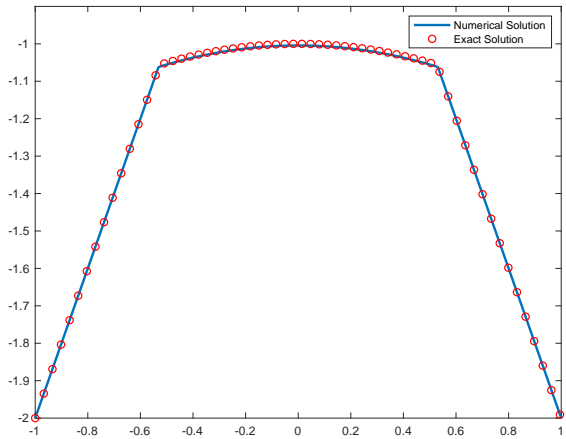}  
\caption{\Allen{Example 4 at time $t = 1.0$ with $m = 3$ using $N = 321$ cells. The solid line is the numerical solution and the dots are the exact solution.} \label{fig:oneDimRiemannPlot}}
\end{center}
\end{figure}
}


\subsection{Examples in Two Dimensions}
\subsubsection{Example 5}
In this example we solve the two-dimensional Burgers' equation
\begin{align*}
\varphi_t + \frac{1}{2}(\varphi_x + \varphi_y)^2 = 0,
\end{align*}
with initial condition \mbox{$\varphi(x,y,0) = -\cos(x+y)$} and periodic boundary conditions on $[0,2\pi]^2$. This equation can be reduced to a one-dimensional equation via the change of variables $z = \frac{x+y}{2}$. That is, 
\begin{align*}
 \frac{\partial u}{\partial z} &= \frac{\partial u}{\partial z} \frac{\partial z}{\partial x}+ \frac{\partial u}{\partial z} \frac{\partial z}{\partial y}\\
 & = \frac{1}{2}\frac{\partial u}{\partial z} + \frac{1}{2}\frac{\partial u}{\partial z}\\
 & =\frac{\partial u}{\partial z}.
\end{align*}
Thus, our equation becomes 
\begin{align*}
\varphi_t + \frac{1}{2}\varphi_z^2 = 0,
\end{align*}
with initial condition \mbox{$\varphi(z,0) = -\cos(2z)$} and periodic boundary conditions on $[0,2\pi]$.

\begin{figure}[htb]
\begin{center}
\includegraphics[width=0.48\textwidth]{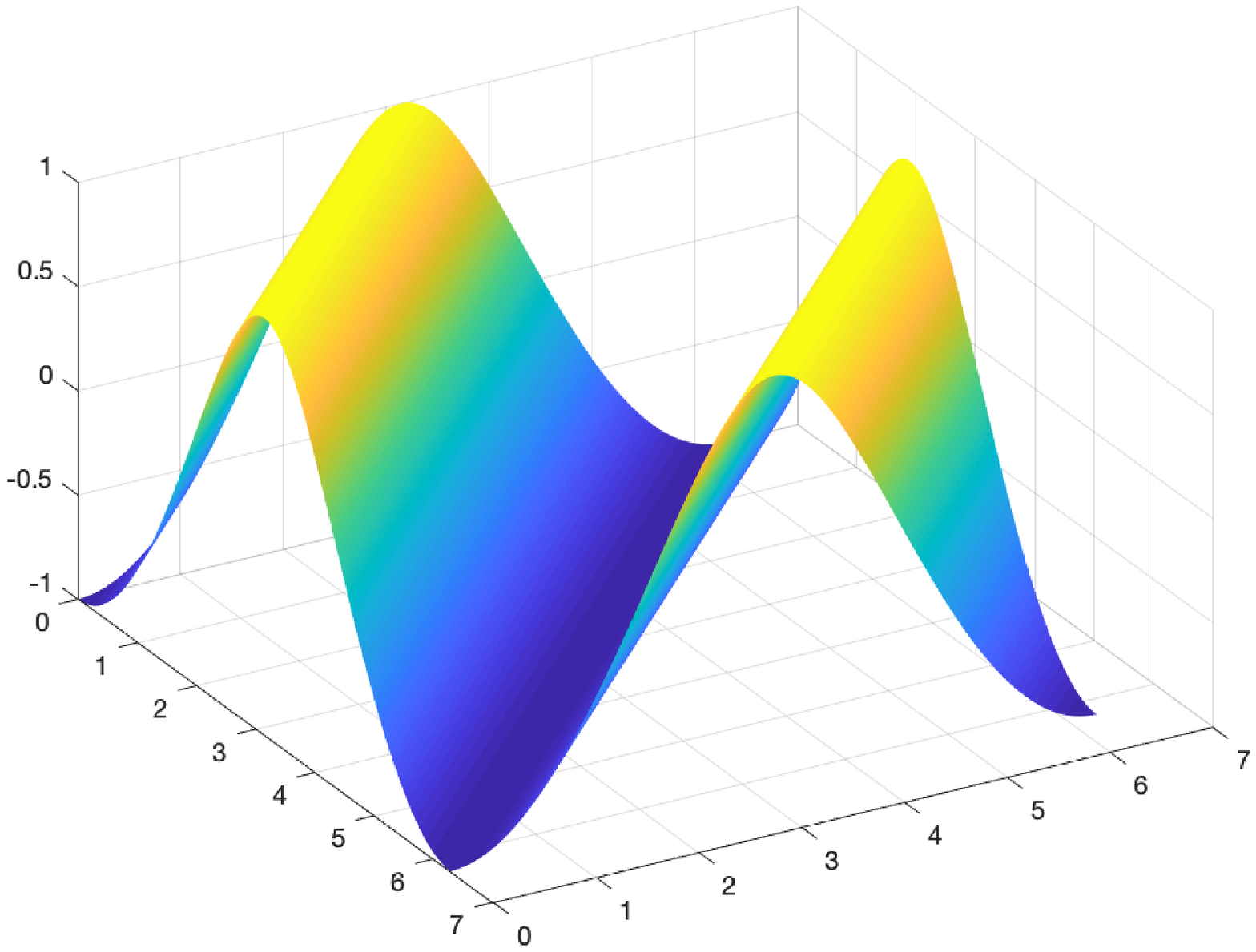}
\includegraphics[width=0.48\textwidth]{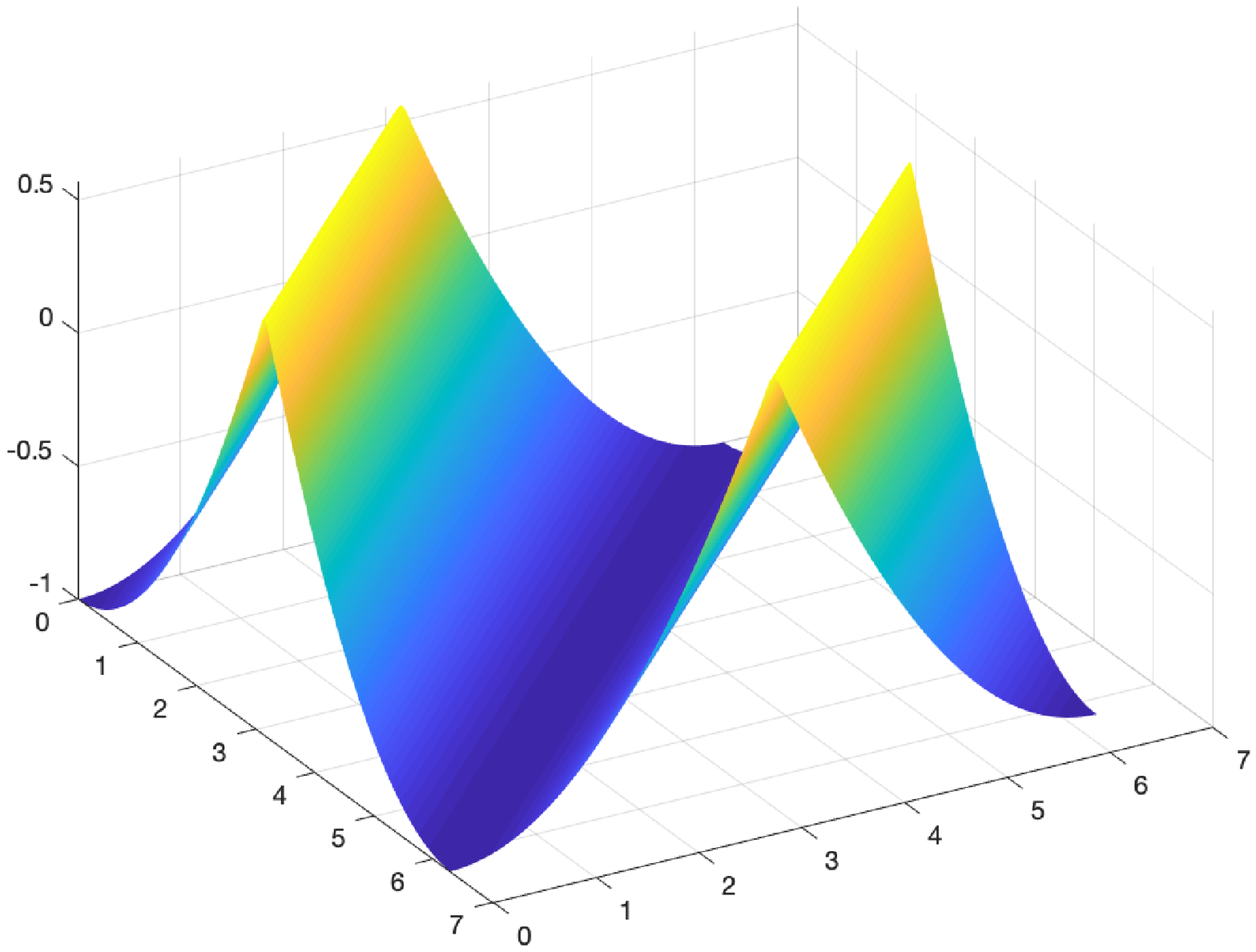} 
\caption{Example 5 on the left is the numerical solution at time $t = 0.1$ approximated using $m = 3$ derivatives on $N = 40$ cells and on the right is the numerical solution at time $t = 0.5$ using the same number of derivatives and cells. \label{fig:2D_burgers}}
\end{center}
\end{figure}

We use the grid with $x_L,y_B = 0$ and $x_R,y_T = 2\pi$ with $n_x,n_y = 10$ cells and refine the grid by a factor of two until we have $n_x,n_y = 80$ cells in order to demonstrate convergence to the viscosity solution. Before the solution develops a kink we demonstrate our method achieves $2m+1$ order accuracy at time $T=0.1$ by giving the $L_1, L_2$ and $L_{\infty}$ norm along with the estimated rates of convergence in Table \ref{tab:2DBurgersSmooth}. We also demonstrate the development of singular features in Figure \ref{fig:2D_burgers}.

We see from the table that we obtain the full order of the method while the solution is smooth. The figure shows us how the method is able to capture the singular features of the solution.
\begin{table}[ht]
 \caption{Errors in Example 5 in the $L_1$, $L_2$ and $L_{\infty}$ norms at time $T = 0.1$ are displayed along with estimated rates of convergence.} 
\begin{center} 
 \begin{tabular}{c c c c c c c } 
  n & $L_1$ error & Conv. Rate & $L_2$ error & Conv. Rate & $L_{\infty}$ error & Conv. Rate \\ 
\hline 
& & & $m = 2$ & & &\\
\hline
 10 & 4.77e-03 & - & 3.78e-04 & - & 8.76e-04 & - \\ 
20 & 1.77e-04 & 4.75 & 1.43e-05 & 4.72 & 3.91e-05 & 4.49 \\ 
40 & 5.72e-06 & 4.95 & 5.77e-07 & 4.63 & 1.33e-06 & 4.88 \\ 
80 & 1.81e-07 & 4.98 & 2.46e-08 & 4.55 & 4.23e-08 & 4.98 \\ 
\hline 
& & & $m = 3$ & & &\\
\hline  
 10 & 1.28e-04 & - & 1.41e-05 & - & 4.29e-05 & - \\ 
20 & 9.70e-07 & 7.05 & 4.90e-08 & 8.17 & 3.62e-07 & 6.89 \\ 
40 & 7.45e-09 & 7.02 & 2.67e-10 & 7.52 & 2.50e-09 & 7.18 \\ 
80 & 5.66e-11 & 7.04 & 2.15e-12 & 6.96 & 1.88e-11 & 7.06 \\ 
\hline
\end{tabular} 
\end{center} 
\label{tab:2DBurgersSmooth}
\end{table} 

\subsubsection{Example 6}
\begin{figure}[htb]
\begin{center}
\includegraphics[width=0.48\textwidth]{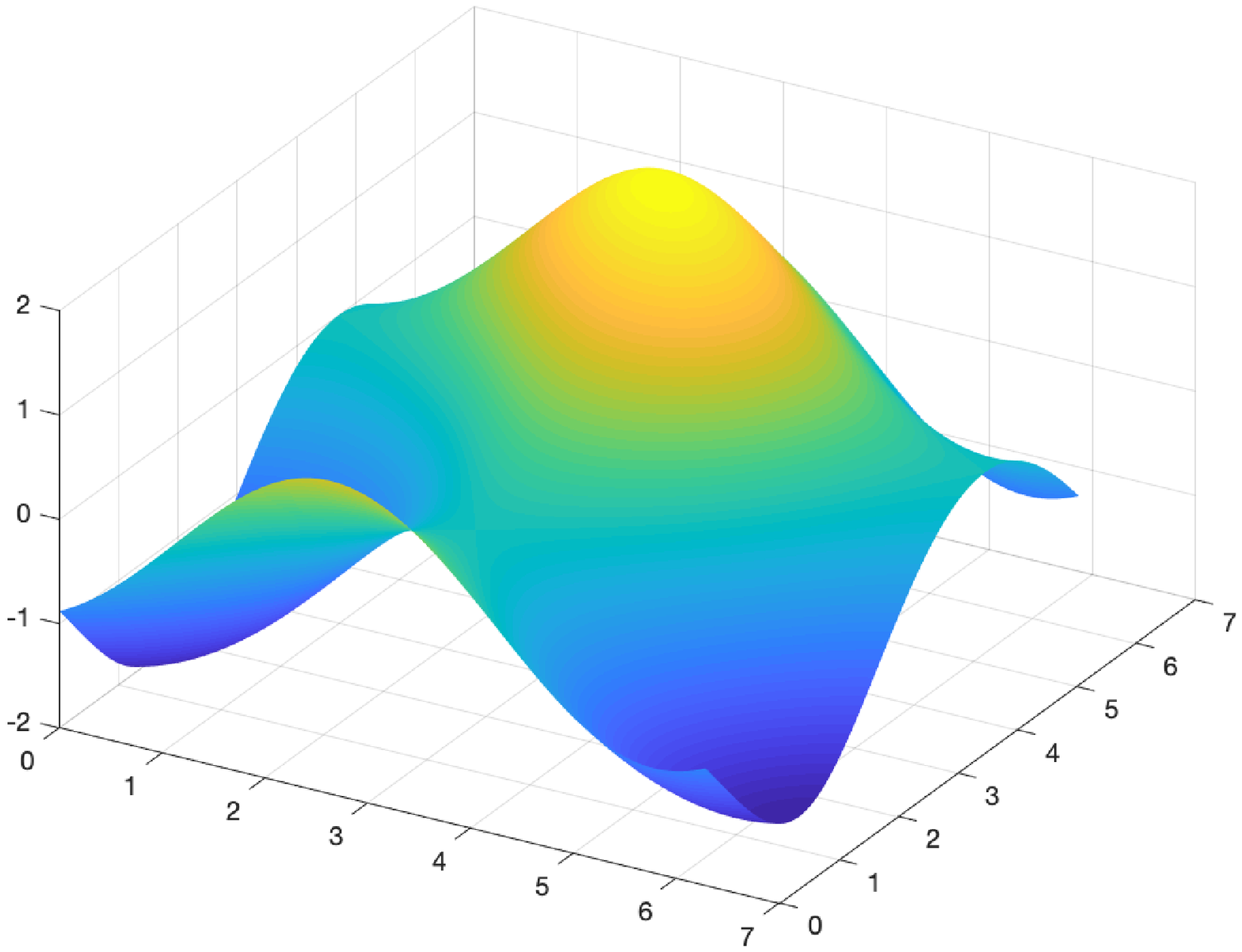} 
\includegraphics[width=0.48\textwidth]{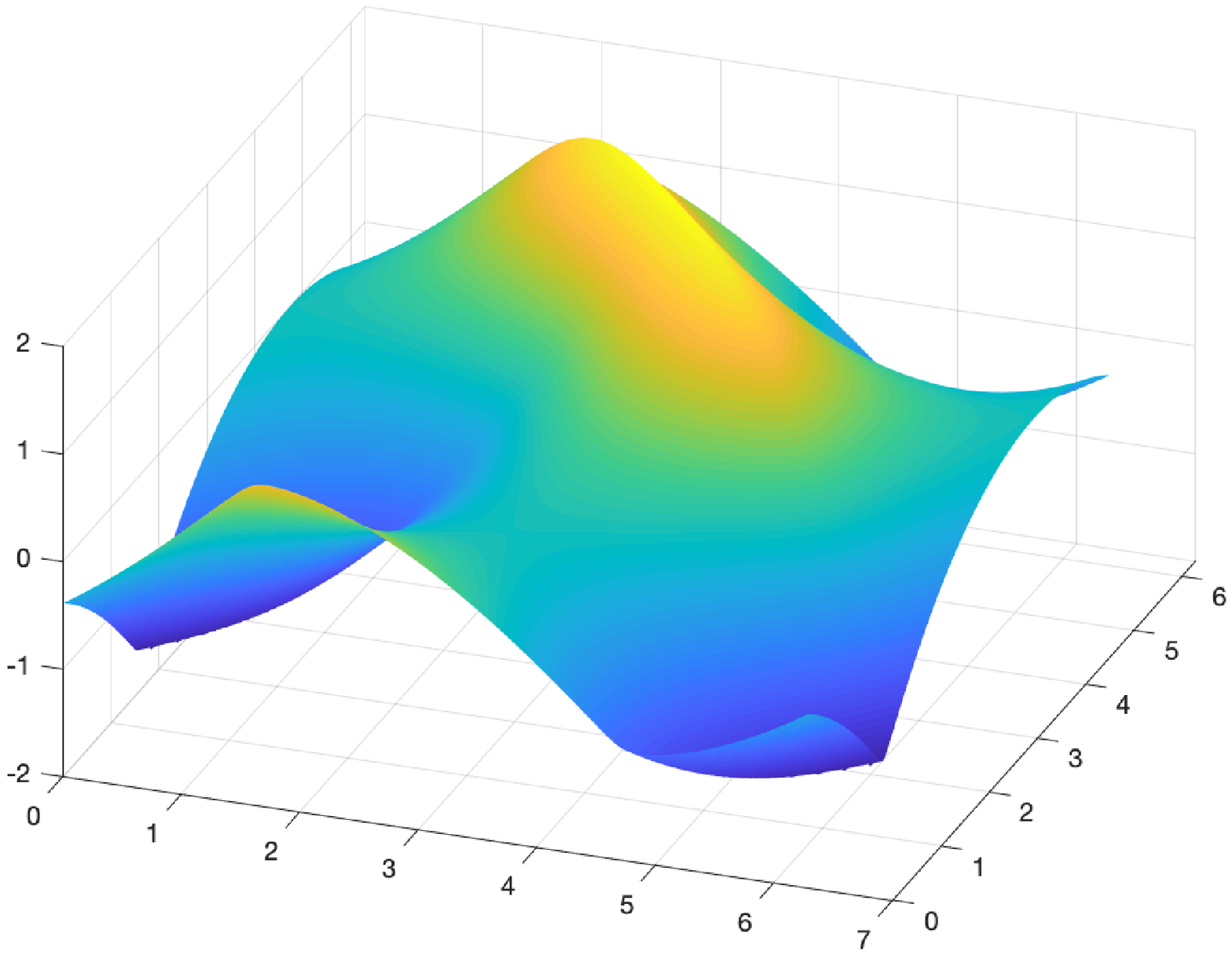} 
\caption{Example 6 on the left is the numerical solution at time $t = 0.5$ approximated using $m = 3$ derivatives on $N = 40$ cells on the right is the numerical solution at time $t = 1.5$ using the same number of derivatives and cells. \label{fig:2D_nonconvex}}
\end{center}
\end{figure}

In this example we solve a two-dimensional nonlinear equation 
\begin{align*} 
\varphi_t + \varphi_x\varphi_y=0,
\end{align*}
with initial condition \mbox{$\varphi(x,y,0) = \sin(x) + \cos(y)$} and periodic boundary conditions on the domain $[-\pi,\pi]^2$.\\\\
This is a genuinely nonlinear problem with a nonconvex Hamiltonian. The viscosity solution is smooth at time $T = 0.5$; we demonstrate $2m+1$ convergence at this time. By $T=1.5$ the viscosity solution develops singular features. We use the grid with $x_L,y_B = -\pi$ and $x_R,y_T = \pi$ with $n_x,n_y = 10$ cells and refine the grid by a factor of two until we have $n_x,n_y = 80$ cells in order to demonstrate convergence to the viscosity solution. Before the solution develops singular features we demonstrate our method achieves $2m+1$ order accuracy at time $T=0.5$ by giving the $L_1, L_2$ and $L_{\infty}$ norm along with the estimated rates of convergence in Table \ref{tab:Example5Smooth}. We also demonstrate the singular features that the solution develops in Figure \ref{fig:2D_nonconvex}
\begin{table}[ht]
 \caption{Errors in Example 6 in the $L_1$, $L_2$ and $L_{\infty}$ norms at time $T = 0.5$ are displayed along with estimated rates of convergence.} 
\begin{center} 
 \begin{tabular}{c c c c c c c } 
\hline 
 n & $L_1$ error & Convergence & $L_2$ error & Convergence & $L_{\infty}$ error & Convergence \\ 
 \hline 
 & & & $m = 2$ & & &\\
 \hline 
 10 & 3.50e-04 & - & 6.78e-05 & - & 2.18e-05 & - \\ 
20 & 1.08e-05 & 5.02 & 2.11e-06 & 5.00 & 6.89e-07 & 4.98 \\ 
40 & 3.35e-07 & 5.01 & 6.58e-08 & 5.00 & 2.14e-08 & 5.01 \\ 
80 & 1.04e-08 & 5.01 & 2.05e-09 & 5.00 & 6.64e-10 & 5.01 \\ 
\hline 
& & & $m = 3$ & & &\\
 \hline 
 10 & 2.62e-07 & - & 5.15e-08 & - & 1.76e-08 & - \\ 
20 & 1.95e-09 & 7.07 & 3.86e-10 & 7.06 & 1.33e-10 & 7.05 \\ 
40 & 1.48e-11 & 7.04 & 2.96e-12 & 7.03 & 1.00e-12 & 7.05 \\ 
80 & 1.17e-13 & 6.98 & 2.52e-14 & 6.87 & 1.51e-14 & 6.05 \\ 
\hline 
\end{tabular} 
\end{center} 
\label{tab:Example5Smooth}
\end{table} 

This example is truly a two-dimensional nonlinear problem and we still see that we obtain the full order of the method while the solution is smooth and our method is able to capture the singular features of the solution.

\Allen{
\subsubsection{Example 7}
In this example we solve a two-dimensional Riemann problem with a nonconvex Hamiltonian
\begin{align*} 
\varphi_t + \sin(\varphi_x + \varphi_y) = 0,
\end{align*}
with initial data $\varphi(x,y,0) = \pi(|y| - |x|)$ on the domain $\Omega = [-1,1] \times [-1,1]$. All of the waves propagate out of the domain and no physical boundary conditions are needed. As a naive implementation of outflow boundary conditions we simply extend sufficiently so that the influence of periodic boundary conditions does not affect the solution during the simulation time. The results of a simulation using $m = 2$ are displayed in Figure \ref{fig:2D_Riemann}.


\begin{figure}[htb]
\begin{center}
\includegraphics[width=0.8\textwidth]{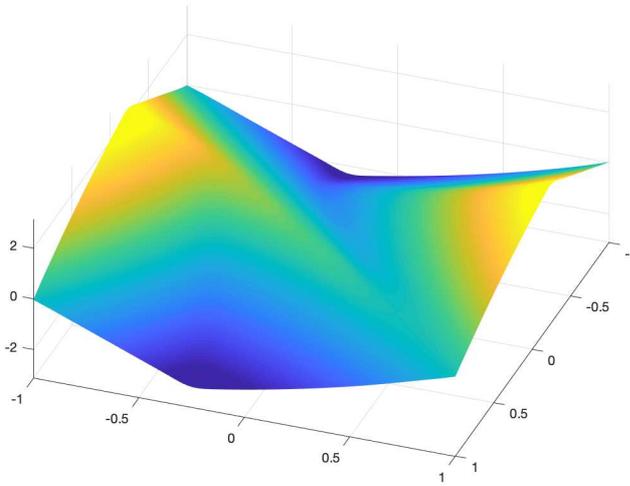}
\caption{\Allen{Example 7. Plotted is the numerical solution at time $t = 1.0$ approximated using $m = 2$ derivatives. The computation use $N = 20$ cells in the physical domain $\Omega = [-1,1] \times [-1,1]$.} \label{fig:2D_Riemann}}
\end{center}
\end{figure} 
}

\Allen{
\subsubsection{Example 8}
In this example we solve a problem related to optimal cost determination 
\begin{align*}
\varphi_t + \sin(y)\varphi_x + (\sin(x)+\text{sign}(\varphi_y))\varphi_y - \frac{1}{2} \sin^2(y) + \cos(x) - 1 = 0,
\end{align*}
with initial data $\varphi(x,y,0) = 0$ and periodic boundary conditions on $\Omega = [-\pi,\pi] \times [-\pi,\pi]$.

The Hamiltonian is not smooth for this example. We are able to capture the viscosity solution well. In Figure \ref{fig:2D_OptimalControl} we display the numerical solution on the left and the optimal control term $\text{sign}(\varphi_y)$ on the right. 

\begin{figure}[htb]
\begin{center}
\includegraphics[width=0.48\textwidth]{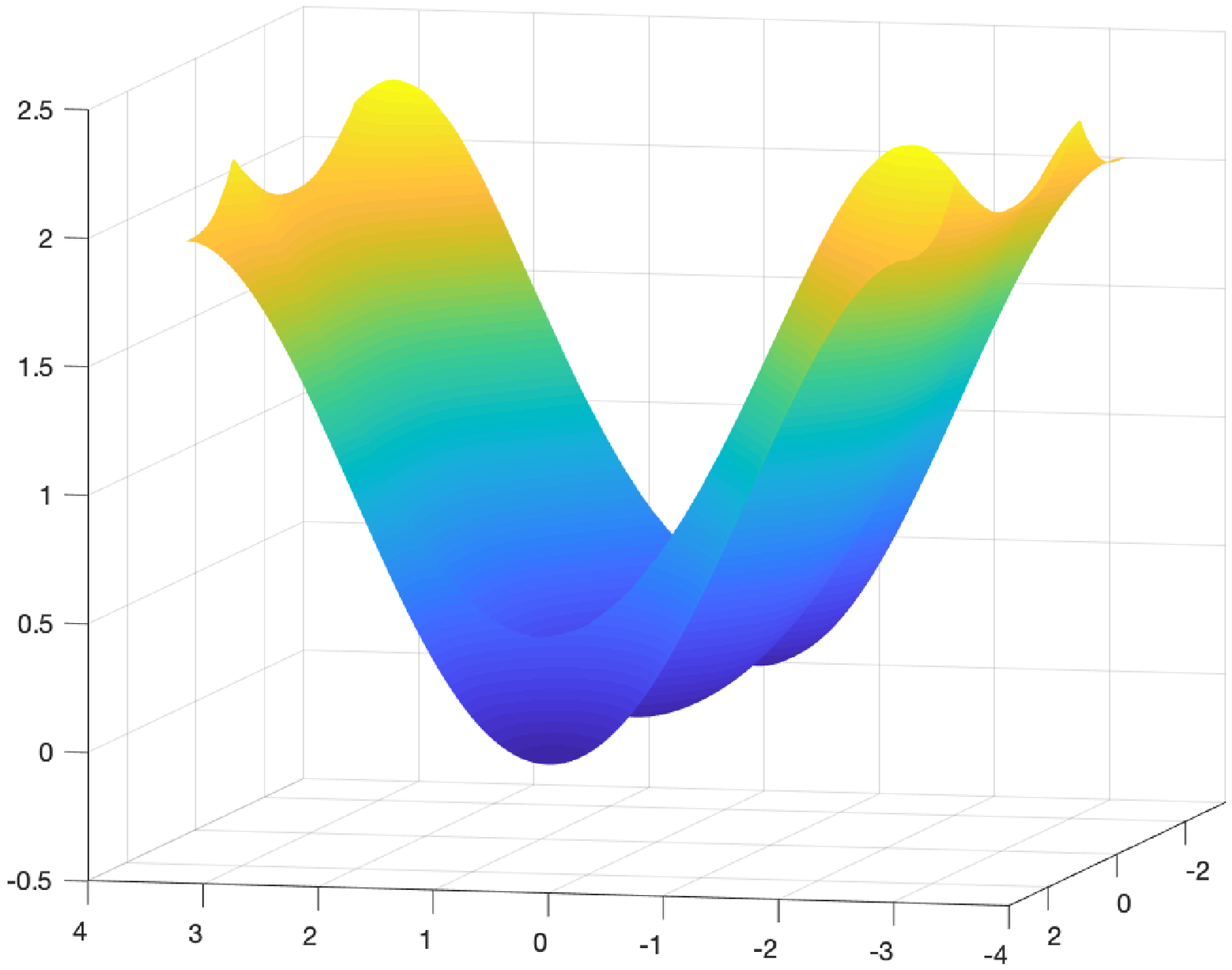}
\includegraphics[width=0.48\textwidth]{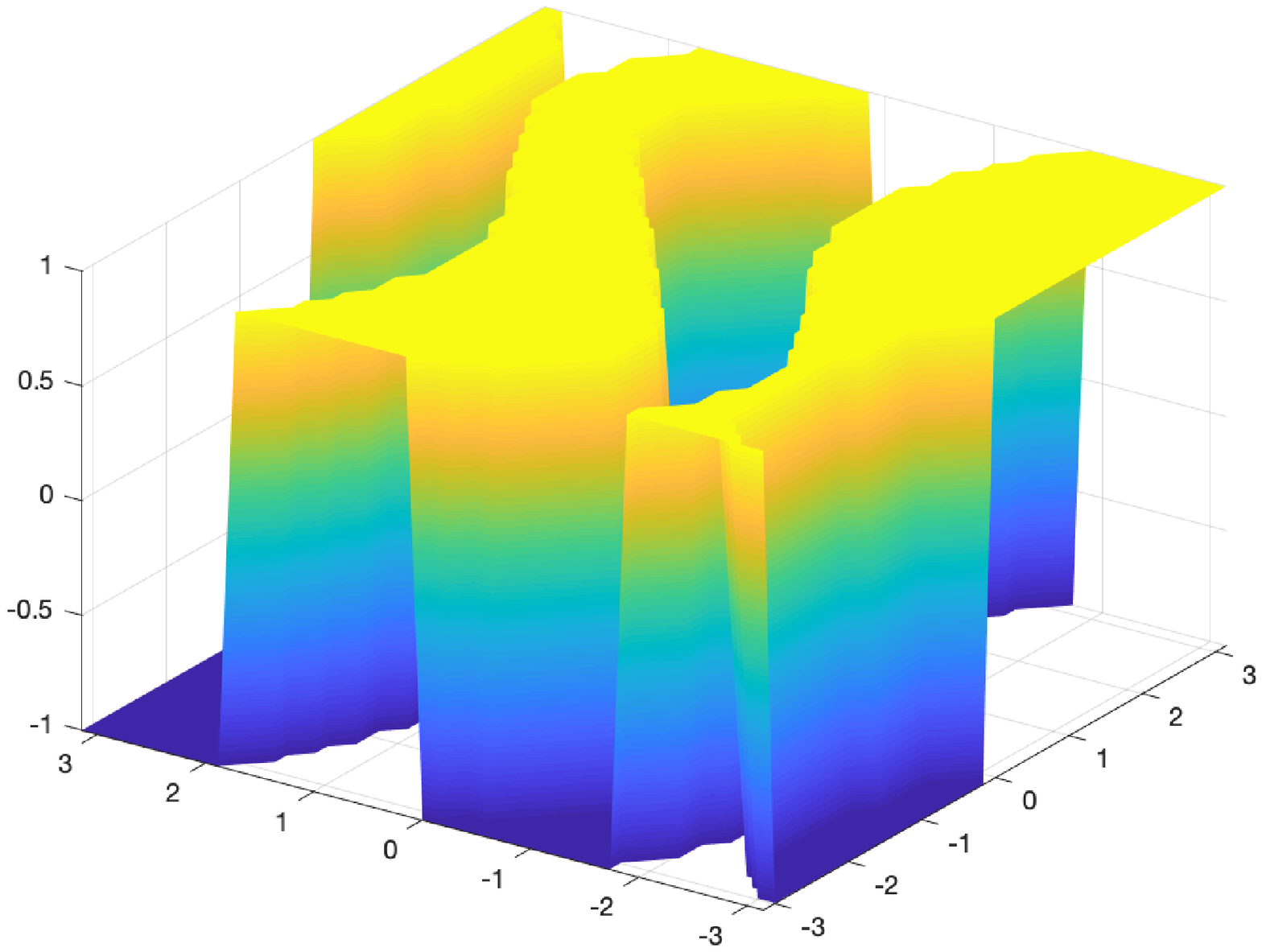}
\caption{\Alldeaa{Example 8 on the left is the numerical solution at time $t = 1.0$ approximated using $m = 2$ derivatives with $N = 40$  cells. On the right is the optimal control term $\text{sign}(\varphi_y)$.  \label{fig:2D_OptimalControl}}}
\end{center}
\end{figure} 
}

\section{Conclusions and Future Work}
Through the coupling of Hermite methods and a discontinuity sensor, our method attains $2m+1$ order of convergence in smooth regions while converging to the viscosity solution when kinks are present. While we were able to achieve the goals of: 1) high order accuracy in smooth regions and 2) sharp resolution of kinks, we believe that there are several ways this method can be improved upon. The Eikonal equation gave us inspiration to develop a flux-conservative Hermite method for HJ equations in order to keep continuity in the derivatives at the element interfaces even when the Hamiltonian is discontinuous. Our method, while effective for a Cartesian grid, can not handle complex geometries. The next step is to deal with different types of boundary conditions and apply Hermite methods to curvilinear coordinate systems. We believe that sensing on a curvilinear coordinate system will be a straightforward generalization since we can map the curvilinear element onto the reference element on the unit square. 

\section{Declarations}
\subsection{Funding}
Allen Alvarez Loya was funded by NSF Grant DGE-1650115. Daniel Appel\"{o} was supported, in part, by NSF Grant DMS-1913076.
\subsection{Conflicts of interest/Competing interests}
On behalf of all authors, the corresponding author states that there is no conflict of interest.
\subsection{Availability of data and material}
All data generated or analyzed during this study are included in this published article.
\subsection{Code availability}
The numerical examples in can be found in the github repository \url{https://github.com/allenalvarezloya/Hermite_HJ}.
\subsection{Authors' contributions}
\FinalEdits{
\section*{Appendix}
\appendix
\section{Solution Methods}
\subsection{Viscosity Solution for a Convex Hamiltonian}
We solve Hamilton-Jacobi problems of the form 
\begin{align*}
u_t + H(Du) &= 0 \quad \text{in } \mathbb{R}^n \times (0,\infty),\\
u &= g \quad \text{on } \mathbb{R}^n \times \lbrace t = 0 \rbrace.
\end{align*}
Here $u:\mathbb{R}^n \times [0,\infty) \to \mathbb{R}$ is the unknown and $Du = D_xu = (u_{x_1}, \dots, u_{x_n})$. We are given the Hamiltonian, $H$, and the initial function $g$.

If the Hamiltonian is convex, then we may use the Lax-Hopf formula given by 
\[
u(x,t) = \min_{y \in \mathbb{R}} \left\lbrace tL\left(\frac{x-y}{t} \right) + g(y) \right \rbrace,
\]
where $g$ is the initial data and $L$ is the Lagrangian. The Lagrangian $L$ and Hamiltonian $H$ are related by the following equations
\begin{align*}
H(p) &= p\cdot \mathbf{v}(p) - L(\mathbf{v}(p)) \quad (p \in \mathbb{R}^n),\\
L(v) &= v\cdot \mathbf{p}(v) - H(\mathbf{p}(v)) \quad (v \in \mathbb{R}^n), 
\end{align*}
where $p = DL(v)$ and $v = DH(p)$.

\subsection{Viscosity Solution for a nonConvex Hamiltonian}
For examples when the Hamiltonian is not convex we can not use the Lax-Hopf formula. We consider the general Hamilton-Jacobi partial differential equation 
\[
u_t + H(Du,x) = 0,
\]
where $Du = D_xu = (u_{x_1},\dots,u_{x_n})$. While the solution is smooth we obtain the solution using the characteristics given by (Details in section 3.2 of \cite{evans2010partial})
\begin{align*}
\dot{\mathbf{p}}(s) &= -D_xH(\mathbf{p}(s),\mathbf{x}(s)),\\
\dot{z}(s) &= D_pH(\mathbf{p}(s),\mathbf{x}(s))\cdot\mathbf{p}(s) - H(\mathbf{p}(s),\mathbf{x}(s)),\\
\dot{\mathbf{x}}(s) &= D_p H(\mathbf{p}(s),\mathbf{x}(s)),
\end{align*} 
for $p(\cdot) = (p^1(\cdot), \dots, p^n(\cdot)), z(\cdot)$, and $\mathbf{x}(\cdot) = (x^1(\cdot), \dots,x^n(\cdot))$.

\subsection{Viscosity Solution for Riemann Problems with Special Initial Data}
We obtain the solution for Riemann problems with special initial data \cite{osher1983riemann}
\[
\varphi_t + f(\varphi_x) + g(\varphi_y) = 0,
\]
with initial data of the form 
\begin{align*}
\varphi(x,y,0) & = (x\sin \theta - y\cos \theta)u^L + (x\cos \theta + y \sin \theta)u^M, \text{ if } x\sin \theta - y \cos \theta \leq 0,\\
\varphi(x,y,0) & = (x\sin \theta - y\cos \theta)u^R + (x\cos \theta + y \sin \theta)u^M, \text{ if } x\sin \theta - y \cos \theta > 0.
\end{align*}
Here $u^R, u^L$ and $\theta$ are fixed constants with $0 \leq \theta < \pi$.

The solution is given by 
\begin{align*}
\varphi(x,y,t) & = (x \sin \theta - y \cos \theta) + t\left(\max \limits_{v \in [u^L,u^R] } \xi v - f(v\sin \theta + u^M\cos \theta ) - g(-v\cos \theta + u^M \sin \theta )\right), \text{ if } u^L < u^R,\\
\varphi(x,y,t) & = (x \sin \theta - y \cos \theta) + t\left(\min \limits_{v \in [u^L,u^R] } \xi v - f(v\sin \theta + u^M\cos \theta ) - g(-v\cos \theta + u^M \sin \theta )\right), \text{ if } u^L > u^R.\\
\end{align*}
where $\xi = \frac{x\sin \theta - y \cos \theta}{t}$.

\subsection{Examples}
\subsubsection{Example 1}
In this example we solve the one-dimensional Burgers' equation 
\begin{equation}
\varphi_t +\frac{1}{2} (\varphi_x)^2=0,
\nonumber
\end{equation}
with initial condition \mbox{$\varphi(x,0)=\sin(x)$} and with periodic boundary conditions \mbox{$\varphi(0,t)=\varphi(2\pi,t)$}.
The Hamiltonian is convex; therefore, we apply the Lax-Hopf formula to obtain the viscosity solution
\[
u(x,t) = \min_{y \in \mathbb{R}} \left \lbrace \frac{(x-y)^2}{2t} + \sin(y) \right \rbrace.
\]

\subsubsection{Example 2}
In this example we solve the one-dimensional Eikonal equation
\begin{align*} 
\varphi_t + |\varphi_x|=0,
\end{align*}
with initial condition \mbox{$\varphi(x,0)=\sin(x)$} and with periodic boundary conditions.
The Hamiltonian is convex; therefore, we apply the Lax-Hopf formula to obtain the viscosity solution. First we compute the Legendre transform of $H$.
\begin{align*}
L(v) &= \sup  \limits_{v \in \mathbb{R}} \lbrace vp - H(p) \rbrace\\
&= \sup  \limits_{v \in \mathbb{R}} \lbrace vp - |p| \rbrace\\
&= \begin{cases} 
      0 & |v| \leq 0 ,\\      
      \infty & \text{ otherwise.}
   \end{cases}\\
&=:  I_{\Omega},
\end{align*}
where in the last step we define $I_{\Omega}$ as the indicator function on $[-1,1]$. The viscosity solution is given by 
\begin{align*}
u(x,t) &= \min_{y \in \mathbb{R}} \left\lbrace tI_\Omega \left(\frac{x-y}{t} \right) + \sin(y) \right \rbrace\\
 &= \min_{y \in [x-t,x+t]} \sin(y).
\end{align*}

\subsubsection{Example 3}
In this example we solve a one-dimensional equation with a nonconvex Hamiltonian
\begin{align*} 
\varphi_t - \cos(\varphi_x + 1) = 0,
\end{align*}
with initial condition $\varphi(x,0)=-\cos(\pi x)$ and periodic boundary conditions $\varphi(-1,t)=\varphi(1,t)$.
The Hamiltonian is not convex so we use method of characteristics to solve while the solution is smooth. The Hamiltonian for this PDE is $H(p) = -\cos(p+1)$. We arrive at the system of ODEs:
\begin{align*}
\dot{p} &= 0,\\
\dot{z} &= \sin(p+1)p + \cos(p+1),\\
\dot{x} &= \sin(p+1).
\end{align*}
Which have solution:
\begin{align*}
p &= p^0,\\
z &= t(\sin(p+1)p + \cos(p+1)) + z^0,\\
x &= t\sin(p+1) + x^0.
\end{align*}
We obtain the solution by the following  
\begin{itemize}
\item Fix a point $x \in \Omega$ and time $t$. We solve for $x^0$ in the equation $x = x(t) = t\sin(p+1) + x^0$.
\item Since $p = u_x$ we can relate $p$ to $x_0$ by $p = p^0 = u_x(x^0,0) = \pi\sin(\pi x^0)$.
\item We can solve for $x^0$ by using solving $x = t\sin(\pi \sin(\pi x^0) + 1) + x^0$ for a fixed x and t.
\item We solve using Newton's method. Once we compute $x^0$ we can compute $p$ and $u(x,t) = z(s)$.
\end{itemize}
\subsubsection{Example 4}
In this example we solve a one-dimensional Riemann problem with a nonconvex Hamiltonian
\begin{align*} 
\varphi_t + \frac{1}{4}(\varphi_x^2 - 1)(\varphi_x^2 - 4) = 0,
\end{align*}
with initial data $\varphi(x,0) = -2|x|$.

We use the solution for Riemann problems with special initial data to obtain our viscosity solution. We have $f(\varphi_x) = \frac{1}{4}(\varphi_x^2 - 1)(\varphi_x^2 - 4), g(\varphi_y) = 0, u^L = 2, u^M = 0, u^R = -2, \theta = \frac{\pi}{2}$ and $\xi = \frac{x}{t}$. The viscosity solution is given by 
\[
\varphi(x,y,t) = t \min \limits_{v \in [u^L,u^R]} \left\lbrace \frac{x}{t} v - \frac{1}{4}(v^2 - 1)(v^2 - 4) \right\rbrace
\] 

\subsubsection{Example 5}
In this example we solve the two-dimensional Burgers' equation
\begin{align*}
\varphi_t + \frac{1}{2}(\varphi_x + \varphi_y)^2 = 0,
\end{align*}
with initial condition \mbox{$\varphi(x,y,0) = -\cos(x+y)$} and periodic boundary conditions on $[0,2\pi]^2$.

This equation can be reduced to a one-dimensional equation via the change of variables $z = \frac{x+y}{2}$. That is, 
\begin{align*}
 \frac{\partial u}{\partial z} &= \frac{\partial u}{\partial z} \frac{\partial z}{\partial x}+ \frac{\partial u}{\partial z} \frac{\partial z}{\partial y}\\
 & = \frac{1}{2}\frac{\partial u}{\partial z} + \frac{1}{2}\frac{\partial u}{\partial z}\\
 & =\frac{\partial u}{\partial z}.
\end{align*}
Thus, our equation becomes 
\begin{align*}
\varphi_t + \frac{1}{2}\varphi_z^2 = 0,
\end{align*}
with initial condition \mbox{$\varphi(z,0) = -\cos(2z)$} and periodic boundary conditions on $[0,2\pi]$.

The Hamiltonian is convex; therefore, we use the Lax-Hopf formula to obtain the viscosity solution 

\[
u(z,t) = \min_{y\in \mathbb{R}} \left \lbrace \frac{(z-y)^2}{2t} - \cos(2y) \right \rbrace,
\]
We then obtain the solution in $x$ and $y$ by setting $x = x$ and $y = 2z - x$.

\subsubsection{Example 6}
In this example we solve a two-dimensional nonlinear equation 
\begin{align*} 
\varphi_t + \varphi_x\varphi_y=0,
\end{align*}
with initial condition \mbox{$\varphi(x,y,0) = \sin(x) + \cos(y)$} and periodic boundary conditions on the domain $[-\pi,\pi]^2$.

We find the solution using the characteristics while the solution is smooth. The characteristics for this problem are given by 
\begin{align*}
\mathbf{p}' &= 0,\\
z' &=  p_1p_2,\\
\mathbf{x}' &= (p_2,p_1).\\
\end{align*}
Which have solution 
\begin{align*}
z &= tp_1p_2 + z_0,\\
\mathbf{x} &= t(p_2,p_1) + \mathbf{x}^0,\\
\mathbf{p} &= \mathbf{p}_0.
\end{align*}
We find the solution by the following 
\begin{itemize}
\item The initial data implies $p_1 = \cos(x^0), p_2 = -\sin(y^0)$. 
\item We can use this to derive a fixed point problem for $\mathbf{x}$ given by $\mathbf{x} = (-t\sin(y^0)+x^0,t\cos(x^0)+y^0)$. 
\item Once we solve the fixed point problem we may solve for  $p_1$ and $p_2$.
\item  We can now compute the viscosity solution $u(x,y,t) = z(s)$.
\end{itemize}
}
\bibliography{HJ_manuscript}
\end{document}